\documentclass[12pt]{article}
\usepackage{latexsym, amssymb}
\usepackage{mathbbol}

\DeclareMathAlphabet\gothic{U}{euf}{m}{n}

\makeatletter
\def\eqnarray{\stepcounter{equation}\let\@currentlabel=\theequation
\global\@eqnswtrue
\tabskip\@centering\let\\=\@eqncr
$$\halign to \displaywidth\bgroup\hfil\global\@eqcnt\z@
  $\displaystyle\tabskip\z@{##}$&\global\@eqcnt\@ne
  \hfil$\displaystyle{{}##{}}$\hfil
  &\global\@eqcnt\tw@ $\displaystyle{##}$\hfil
  \tabskip\@centering&\llap{##}\tabskip\z@\cr}

\def\endeqnarray{\@@eqncr\egroup
      \global\advance\c@equation\m@ne$$\global\@ignoretrue}

\def\@yeqncr{\@ifnextchar [{\@xeqncr}{\@xeqncr[5pt]}}
\makeatother

\begin{document}
\bibliographystyle{tom}

\newtheorem{lemma}{Lemma}[section]
\newtheorem{thm}[lemma]{Theorem}
\newtheorem{cor}[lemma]{Corollary}
\newtheorem{voorb}[lemma]{Example}
\newtheorem{rem}[lemma]{Remark}
\newtheorem{prop}[lemma]{Proposition}
\newtheorem{stat}[lemma]{{\hspace{-5pt}}}

\newenvironment{remarkn}{\begin{rem} \rm}{\end{rem}}
\newenvironment{exam}{\begin{voorb} \rm}{\end{voorb}}

\newcommand{\gota}{\gothic{a}}
\newcommand{\gotb}{\gothic{b}}
\newcommand{\gotc}{\gothic{c}}
\newcommand{\gote}{\gothic{e}}
\newcommand{\gotf}{\gothic{f}}
\newcommand{\gotg}{\gothic{g}}
\newcommand{\gothh}{\gothic{h}}
\newcommand{\gotk}{\gothic{k}}
\newcommand{\gotm}{\gothic{m}}
\newcommand{\gotn}{\gothic{n}}
\newcommand{\gotp}{\gothic{p}}
\newcommand{\gotq}{\gothic{q}}
\newcommand{\gotr}{\gothic{r}}
\newcommand{\gots}{\gothic{s}}
\newcommand{\gotu}{\gothic{u}}
\newcommand{\gotv}{\gothic{v}}
\newcommand{\gotw}{\gothic{w}}
\newcommand{\gotz}{\gothic{z}}
\newcommand{\gotA}{\gothic{A}}
\newcommand{\gotB}{\gothic{B}}
\newcommand{\gotG}{\gothic{G}}
\newcommand{\gotL}{\gothic{L}}
\newcommand{\gotS}{\gothic{S}}
\newcommand{\gotT}{\gothic{T}}

\newcounter{teller}
\renewcommand{\theteller}{\Roman{teller}}
\newenvironment{tabel}{\begin{list}%
{\rm \bf \Roman{teller}.\hfill}{\usecounter{teller} \leftmargin=1.1cm
\labelwidth=1.1cm \labelsep=0cm \parsep=0cm}
                      }{\end{list}}

\newcounter{tellerr}
\renewcommand{\thetellerr}{(\roman{tellerr})}
\newenvironment{subtabel}{\begin{list}%
{\rm  (\roman{tellerr})\hfill}{\usecounter{tellerr} \leftmargin=1.1cm
\labelwidth=1.1cm \labelsep=0cm \parsep=0cm}
                         }{\end{list}}

\newcommand{\Ni}{{\bf N}}
\newcommand{\Ri}{{\bf R}}
\newcommand{\Ci}{{\bf C}}
\newcommand{\Ti}{{\bf T}}
\newcommand{\Zi}{{\bf Z}}
\newcommand{\Fi}{{\bf F}}

\newcommand{\proof}{\mbox{\bf Proof} \hspace{5pt}} 
\newcommand{\remark}{\mbox{\bf Remark} \hspace{5pt}}
\newcommand{\ruimte}{\vskip10.0pt plus 4.0pt minus 6.0pt}

\newcommand{\simh}{{\stackrel{{\rm cap}}{\sim}}}
\newcommand{\ad}{{\mathop{\rm ad}}}
\newcommand{\Ad}{{\mathop{\rm Ad}}}
\newcommand{\Aut}{\mathop{\rm Aut}}
\newcommand{\arccot}{\mathop{\rm arccot}}
\newcommand{\capp}{{\mathop{\rm cap}}}
\newcommand{\rcapp}{{\mathop{\rm rcap}}}
\newcommand{\diam}{\mathop{\rm diam}}
\newcommand{\divv}{\mathop{\rm div}}
\newcommand{\codim}{\mathop{\rm codim}}
\newcommand{\RRe}{\mathop{\rm Re}}
\newcommand{\IIm}{\mathop{\rm Im}}
\newcommand{\Tr}{{\mathop{\rm Tr}}}
\newcommand{\Vol}{{\mathop{\rm Vol}}}
\newcommand{\card}{{\mathop{\rm card}}}
\newcommand{\supp}{\mathop{\rm supp}}
\newcommand{\sgn}{\mathop{\rm sgn}}
\newcommand{\essinf}{\mathop{\rm ess\,inf}}
\newcommand{\esssup}{\mathop{\rm ess\,sup}}
\newcommand{\Int}{\mathop{\rm Int}}
\newcommand{\Leibniz}{\mathop{\rm Leibniz}}
\newcommand{\lcm}{\mathop{\rm lcm}}
\newcommand{\loc}{{\rm loc}}

\newcommand{\mod}{\mathop{\rm mod}}
\newcommand{\spann}{\mathop{\rm span}}
\newcommand{\one}{\mathbb{1}}

\hyphenation{groups}
\hyphenation{unitary}

\newcommand{\tfrac}[2]{{\textstyle \frac{#1}{#2}}}

\newcommand{\cb}{{\cal B}}
\newcommand{\cc}{{\cal C}}
\newcommand{\cd}{{\cal D}}
\newcommand{\ce}{{\cal E}}
\newcommand{\cf}{{\cal F}}
\newcommand{\ch}{{\cal H}}
\newcommand{\ci}{{\cal I}}
\newcommand{\ck}{{\cal K}}
\newcommand{\cl}{{\cal L}}
\newcommand{\cm}{{\cal M}}
\newcommand{\co}{{\cal O}}
\newcommand{\cs}{{\cal S}}
\newcommand{\ct}{{\cal T}}
\newcommand{\cx}{{\cal X}}
\newcommand{\cy}{{\cal Y}}
\newcommand{\cz}{{\cal Z}}

\newlength{\hightcharacter}
\newlength{\widthcharacter}
\newcommand{\covsup}[1]{\settowidth{\widthcharacter}{$#1$}\addtolength{\widthcharacter}{-0.15em}\settoheight{\hightcharacter}{$#1$}\addtolength{\hightcharacter}{0.1ex}#1\raisebox{\hightcharacter}[0pt][0pt]{\makebox[0pt]{\hspace{-\widthcharacter}$\scriptstyle\circ$}}}
\newcommand{\cov}[1]{\settowidth{\widthcharacter}{$#1$}\addtolength{\widthcharacter}{-0.15em}\settoheight{\hightcharacter}{$#1$}\addtolength{\hightcharacter}{0.1ex}#1\raisebox{\hightcharacter}{\makebox[0pt]{\hspace{-\widthcharacter}$\scriptstyle\circ$}}}
\newcommand{\scov}[1]{\settowidth{\widthcharacter}{$#1$}\addtolength{\widthcharacter}{-0.15em}\settoheight{\hightcharacter}{$#1$}\addtolength{\hightcharacter}{0.1ex}#1\raisebox{0.7\hightcharacter}{\makebox[0pt]{\hspace{-\widthcharacter}$\scriptstyle\circ$}}}

\thispagestyle{empty}

\vspace*{1cm}
\begin{center}
{\Large\bf Diffusion determines the manifold} \\[5mm]

\large W. Arendt$^1$, M. Biegert$^1$ and A.F.M. ter Elst$^2$

\end{center}

\vspace{5mm}

\begin{center}
{\bf Abstract}
\end{center}

\begin{list}{}{\leftmargin=1.8cm \rightmargin=1.8cm \listparindent=10mm 
   \parsep=0pt}
\item
We prove under a weak smoothness condition that 
two Riemannian manifold are isomorphic if and only 
there exists an order isomorphism which intertwines with the 
Dirichlet type heat semigroups on the manifolds.

\end{list}

\vspace{4cm}
\noindent
March 2008

\vspace{5mm}

\vspace{15mm}

\noindent
{\bf Home institutions:}    \\[3mm]
\begin{tabular}{@{}cl@{\hspace{10mm}}cl}
1. & Abteilung Angewandte Analysis & 
2. & Department of Mathematics  \\
& Universit\"at Ulm & 
  & University of Auckland  \\
& Helmholtzstr.\ 18  &
  & Private bag 92019 \\
& 89069 Ulm  &
   & Auckland  \\ 
& Germany  &
  & New Zealand \\[8mm]
\end{tabular}

\newpage
\setcounter{page}{1}

\section{Introduction} \label{Sdrum1}

A fundamental problem raised in Kac's famous article \cite{Kac}
`Can one hear the shape of a drum' is whether two isospectral manifolds
are isomorphic.
The answer is negative in general.
Milnor gave a counter example for compact Riemannian manifolds~\cite{Milnor}.
In the Euclidean case the first example was given in dimension~$4$ by 
Urakawa~\cite{Ura}.
Then Gordon--Webb--Wolpert 
\cite{GWW2} constructed two polygons in $\Ri^2$ which are isospectral 
but not isomorphic.
Moreover, \cite{GordonW} constructed two isospectral convex open sets
in $\Ri^4$ which are isospectral but not isomorphic.
Kac's question in the strict sense, namely whether two isospectral 
bounded open sets in $\Ri^2$ with $C^\infty$-boundary are isometric,
is still to be open.
But there are recent positive results by Zelditch \cite{Zel}
for open sets in $\Ri^2$
with analytic boundary verifying some symmetry conditions.

To say that the two manifolds are isospectral means by definition
that the corresponding Dirichlet Laplacians have the same 
eigenvalues counted with multiplicity.
This, in turn, can be reformulated by saying that  there exists a unitary
operator $U$ intertwining the two heat semigroups.
The heat semigroups are positive, i.e.\ positive initial values
lead to positive solutions.
These positive solutions describe the heat diffusion on the manifold.
Thus, if instead of a unitary operator, we consider an order isomorphism
$U$ (i.e.\ $U$ is linear, bijective and $U \varphi \geq 0$ if and only if $\varphi \geq 0$)
on $L_2$, then to say say that $U$ intertwines the heat semigroups
means that $U$ maps the positive solutions to positive solutions.
It was shown in \cite{Are3} that in the Euclidean case,
i.e.\ if we consider open connected sets in $\Ri^d$, then 
these sets are necessarily congruent as soon as such an intertwining 
order isomorphism exists. 
This may be rephrased by saying that 
diffusion determines the body.
The aim of this paper is to extend this result to manifolds.

There are several notable new features coming into play in the 
non-Euclidean case.
First of all, in \cite{Are3} a precise regularity  condition has been established
under which the result is valid.
The open sets have to be regular in capacity
(this means loosely speaking that they do not have holes of
capacity~$0$).
Some effort is made in this paper to extend this notion to manifolds,
which is not possible in an immediate way.

The problem addressed in this paper is partially motivated by 
work of Arveson \cite{Arv2} \cite{Arv1}, who 
introduces differential structures in operator algebras.
Our results imply uniqueness of these differential structures, the case 
of compact Riemannian manifolds being of particular interest.

Not all results in the Euclidean case carry over to Riemannian manifolds.
We give an example, Example~\ref{xdrum501}, of a non-zero lattice homomorphism 
which intertwine the heat semigroups, but which is not an isomorphism,
in contrast to the Euclidean case \cite{Are4}, Theorem~2.1.


\smallskip

Let $(M,g)$ be a Riemannian manifold of dimension $d$.
We always assume that a Riemannian manifold is $\sigma$-compact.
Then $M$ has a natural Radon measure denoted by $|\cdot|$. 
Set 
\[
H^1_\loc(M)
= \{ \varphi \in L_{2,\loc}(M) : \varphi \circ x^{-1} \in H^1_\loc(x(V))
       \mbox{ for every chart } (V,x) \} 
\;\;\; .  \]
If $\varphi \in H^1_\loc(M)$ and $(V,x)$ is a chart on $M$ then set 
$\frac{\partial}{\partial x^i} \varphi = (D_i (\varphi \circ x^{-1})) \circ x \in L_{2,\loc}(V)$,
where $D_i$ denotes the partial derivative in $\Ri^d$.
Moreover, for all $\varphi,\psi \in H^1_\loc(M)$
there exists a unique element $\nabla \varphi \cdot \nabla \psi \in L_{1,\loc}(M)$
such that 
\[
\nabla \varphi \cdot \nabla \psi \Big|_V
= \sum_{i,j=1}^d g^{ij} \Big( \frac{\partial}{\partial x^i} \varphi \Big)
                               \Big( \frac{\partial}{\partial x^j} \psi \Big)
\]
for every chart $(V,x)$ on $M$.
We let $|\nabla \varphi| = (\nabla \varphi \cdot \nabla \varphi)^{1/2}$.
Let $H^1(M)$ be the Hilbert space of all $\varphi \in H^1_\loc(M)$ such that 
both $\varphi, |\nabla \varphi| \in L_2(M)$, with norm
$\varphi \mapsto ( \|\varphi\|_2^2 + \| \, |\nabla \varphi| \, \|_2^2 )^{1/2}$.
Moreover, let $H^1_0(M)$ be the closure of $C_c^\infty(M)$ in $H^1(M)$.
Define the bilinear form $a \colon H^1_0(M) \times H^1_0(M) \to \Ri$ by
$a(\psi,\varphi) = \int \nabla \psi \cdot \nabla \varphi$.
Then $a$ is  closed and positive.
The {\bf Dirichlet Laplace--Beltrami operator} $\Delta$
on $M$ is the associated self-adjoint operator.
If $(V,x)$ is a chart on $M$ then 
\[
\Delta \, \varphi
= - \sum_{i,j=1}^d \frac{1}{\sqrt{g}} \, \frac{\partial}{\partial x^i} \, 
           g^{ij} \, \sqrt{g} \, \frac{\partial}{\partial x^j} \, \varphi
\]
for all $\varphi \in C_c^\infty(V)$.

If $(M_1,g_1)$ and $(M_2,g_2)$ are two Riemannian manifolds 
then a map $\tau \colon M_1 \to M_2$ is called an {\bf isometry}
if it is a $C^\infty$-diffeomorphism and 
\[
g_2|_{\tau(p)}(\tau_*(v), \tau_*(w)) = g_1|_p(v,w)
\]
for all $p \in M_1$ and $v,w \in T_p M$.
A map $\tau \colon M_1 \to M_2$ is called a {\bf local isometry}
if for all $p \in M_1$ there exists an open neighbourhood 
$V$ of $p$ such that the restriction $\tau|_V \colon V \to \tau(V)$
is an isometry.
The Riemannian manifolds $(M_1,g_1)$ and $(M_2,g_2)$
are called {\bf isomorphic} if there exists an isometry 
from $M_1$ onto $M_2$.
If $\tau \colon M_1 \to M_2$ is an isometry
and $p \in [1,\infty)$ then 
$\varphi \circ \tau \in L_p(M_1)$ and
\begin{equation}
\|\varphi \circ \tau\|_{L_p(M_1)} = \|\varphi\|_{L_p(M_2)}
\label{eSdrum1;1}
\end{equation}
for all $\varphi \in L_p(M_2)$.
In particular, the map $\varphi \mapsto \varphi \circ \tau$ is a 
unitary map from $L_2(M_2)$ onto $L_2(M_1)$ and a unitary map
from $H^1_0(M_2)$ onto $H^1_0(M_1)$.
Moreover, if $\varphi \in L_2(M_2)$ then 
$\varphi \in D(\Delta_2)$ if and only if $\varphi \circ \tau \in D(\Delta_1)$
and $\Delta_1(\varphi \circ \tau) = (\Delta_2 \varphi) \circ \tau$, 
where $\Delta_j$ is the Dirichlet Laplace--Beltrami operator 
on $M_j$ for all $j \in \{ 1,2 \} $.

A linear operator $U \colon E \to F$ between two Riesz spaces is said to be a 
{\bf lattice homomorphism} if 
\[
U (\varphi \wedge \psi) = (U \varphi) \wedge (U \psi)
\]
for all $\varphi,\psi \in E$.
For alternative equivalent definitions see \cite{AB} Theorem~7.2.
Here in this paper in most cases the spaces $E$ and $F$ will be $L_p$-spaces
and then $(\varphi \wedge \psi)(x) = \min \{ \varphi(x),\psi(x) \} $ a.e.
Each lattice homomorphism $U$ is positive, i.e.\ $\varphi \geq 0$
implies $U \varphi \geq 0$.
An {\bf order isomorphism} $U \colon E \to F$ is a bijective 
mapping  such that $U \varphi \geq 0$ if and only if $\varphi \geq 0$.
Equivalently, $U$ is an order isomorphism
if and only if $U$ is a bijective lattice homomorphism.
Then also $U^{-1}$ is an order isomorphism.
Recall that also  each positive operator between $L_p$-spaces is continuous 
by \cite{AB} Theorem~12.3.


The main theorem of this paper is the following.
It is valid under some regularity assumptions on the manifolds,
namely regularity in capacity, which is optimal for this purpose and 
which we will explain below.

\begin{thm} \label{tdrum101}
Let $(M_1,g_1)$ and $(M_2,g_2)$ be two connected Riemannian manifolds which are both regular in capacity.
Let $p \in [1,\infty)$.
For all $j \in \{ 1,2 \} $ let $\Delta_j$ be the Dirichlet Laplace--Beltrami operator 
on $M_j$ and let $S^{(j)}$ be the associated semigroup on $L_p(M_j)$.
Then the following two conditions are equivalent.
\begin{tabel}
\item \label{tdrum101-1}
$(M_1,g_1)$ and $(M_2,g_2)$ are isomorphic.
\item \label{tdrum101-3}
There exists a lattice homomorphism $U \colon L_p(M_1) \to L_p(M_2)$ such that 
$U L_p(M_1)$ is dense in $L_p(M_2)$ and 
\[
U S^{(1)}_t = S^{(2)}_t U
\]
for all $t > 0$.
\end{tabel}
Moreover, if $U$ is a lattice homomorphism as in Condition~{\rm \ref{tdrum101-3}} then 
$U$ is an order isomorphism and
there exist $c>0$ and an isometry $\tau \colon M_2 \to M_1$ such that 
$U \varphi = c \, \varphi \circ \tau$ for all $\varphi \in L_p(M_1)$.
\end{thm}

It turns out that 
all complete connected Riemannian manifolds, and in particular all compact
connected Riemannian manifolds, are regular in capacity.
Therefore one immediately has the following corollary.

\begin{cor} \label{cdrum102}
Let $(M_1,g_1)$ and $(M_2,g_2)$ be two complete connected Riemannian manifolds.
Let $p \in [1,\infty)$.
For all $j \in \{ 1,2 \} $ let $\Delta_j$ be the Dirichlet Laplace--Beltrami operator 
on $M_j$ and let $S^{(j)}$ be the associated semigroup on $L_p(M_j)$.
Then the following two conditions are equivalent.
\begin{tabel}
\item \label{cdrum102-1}
$(M_1,g_1)$ and $(M_2,g_2)$ are isomorphic.
\item \label{cdrum102-2}
There exists an order isomorphism $U \colon L_p(M_1) \to L_p(M_2)$ such that 
\[
U S^{(1)}_t = S^{(2)}_t U
\]
for all $t > 0$.
\end{tabel}
\end{cor}

Now we explain the notion of regularity in capacity for Riemannian manifolds.
The {\bf capacity} of a subset $A$ of $M$ is given by
\[
\capp_M(A)
= \capp(A)
= \inf \{ \|\varphi\|_{H^1(M)}^2 : \varphi \in H^1(M) \mbox{ and } \varphi \geq 1
                                 \mbox{ on a neighbourhood of } A \} 
\;\;\; .  \]
An open subset $\Omega$ of $\Ri^d$ is called {\bf regular in capacity}
\cite{Are3} if $\capp_{\Ri^d} (B(x\,;r) \setminus \Omega) > 0$
for all $x \in \partial \Omega$ and $r > 0$, where $B(x\,;r)$ is the 
Euclidean ball. 
Biegert and Warma gave several characterizations for regular in capacity.
In particular, an open subset $\Omega$ of $\Ri^d$ is regular in capacity 
if and only if every $\varphi \in H^1_0(\Omega) \cap C(\overline{\Omega})$ is 
zero everywhere on $\partial \Omega$ (\cite{BW} Theorem~3.2).
Since $\Ri^d$ is locally compact it then follows that 
an open subset $\Omega$ of $\Ri^d$ is regular in capacity 
if and only if every $\varphi \in H^1_0(\Omega) \cap C_0(\overline{\Omega})$ is 
zero everywhere on $\partial \Omega$.
This characterization allows an extension to general connected Riemannian manifolds.
There is a natural distance $d_M$ on a connected Riemannian manifold~$M$.
We denote by $B_M(p\,;r) = B(p\,;r)$ the associated balls.
Let $\widetilde M$ denote the (metric) completion of $M$ with respect to this distance.
Set 
\[ 
\partial M=\widetilde M\setminus M 
\;\;\; .  \]
We say that a connected Riemannian manifold $M$ is {\bf regular in capacity} if
$\varphi(p) = 0$ for all 
$\varphi \in C_0(\widetilde M) \cap H^1_0(M) 
      = \{ \varphi \in C_0(\widetilde M) : \varphi|_M \in H^1_0(M) \} $ and $p \in \partial M$.
Here $C_0(\widetilde M)$ is the closure of the space $C_c(\widetilde M)$
of all continuous functions with compact support, with respect to the 
supremum norm in the space of all bounded continuous functions on~$\widetilde M$.
Clearly every complete connected Riemannian manifold is regular in capacity.

In the Euclidean case regularity in capacity is a very mild  condition
on the boundary of an open subset.
If $\Omega \subset \Ri^d$ is open and bounded then it is regular in capacity
if it is Dirichlet regular. 
The Lebesgue cusp is regular in capacity, but not Dirichlet regular
(see \cite{AD} Section~7).

If $M_1$ and $M_2$ are two isomorphic connected Riemannian manifolds and
$\tau \colon M_1 \to M_2$ is an isometry 
then $\tau$ is {\bf distance preserving}, i.e.\
$d_{M_2}(\tau(p)\,;\tau(q)) = d_{M_1}(p\,;q)$ for all $p,q \in M_1$.
Moreover, if $M_1$ is 
regular in capacity, then also $M_2$ is regular in capacity.

Now we can explain why regularity in capacity is the minimal regularity condition
in our context.
Let $M$ be a connected Riemannian manifold which is 
complete (or more general, regular in capacity).
Let $\emptyset \neq N \subset M$ be a closed subset of capacity zero.
Then $\Omega := M \setminus N$ is again a connected Riemannian manifold
(see Theorem~\ref{tdrum302})
The injection $\tau \colon \Omega \to M$ defines an isometry which is not 
surjective.
The unitary operator $U \colon L_2(M) \to L_2(\Omega)$ given by $U \varphi = \varphi \circ \tau$
is an order isomorphism intertwining the two heat semigroups
even though $\Omega$ and $M$ are not isomorphic.
It follows from Theorem~\ref{tdrum101} that $\Omega$ is not regular in capacity.

The paper is organized as follows.
In the next section we give a sufficient condition to ensure that the 
distance on a subriemannian manifold equals the induced distance.
In Section~\ref{Sdrum5} we show that $M_1$ and $M_2$ are 
isometric if they have sufficiently big isometric open subsets.
In Section~\ref{Sdrum2} we prove Theorem~\ref{tdrum101}.
Finally, in Section~\ref{Sdrum4} we give several characterizations of 
regularity in capacity.

\subsection*{Acknowledgements}
The first and second named authors would like to thank for the great 
hospitality and generosity 
during their stay at the University of Auckland.
The third named author is most grateful for the hospitality extended
to him during a most enjoyable and fruitful stay at the University of Ulm.

\section{Distances} \label{Sdrum3}

If $N$ is a connected open subset of a connected Riemannian manifold $M$
then $d_M(p\,;q) \leq d_N(p\,;q)$ for all $p,q \in N$, where
$d_M$ and $d_N$ are the natural distances on $M$ and $N$.
Even if $|M \setminus N| = 0$, then it is easy to construct examples 
such that the induced distance from $d_M$ on $N$ differs from the 
distance $d_N$.
We next show that the condition $\capp_M(M \setminus N) = 0$ suffices to have equality.

\begin{thm} \label{tdrum302}
Let $N$ be an open subset of a connected Riemannian manifold $M$
and suppose  that $\capp_M(M \setminus N) = 0$.
Then $N$ is connected and $d_M(p\,;q) = d_N(p\,;q)$ for all $p,q \in N$.
\end{thm}

The proof involves an alternative description of the distances.
First we need $L_\infty$-versions of $H^1$ and $\nabla \varphi$.
Let $M$ be a Riemannian manifold.
Set 
\[
W^{1,\infty}_\loc(M) 
= \{ \varphi \in L_{\infty,\loc}(M) : \varphi \circ x^{-1} \in W^{1,\infty}_\loc(x(V))
       \mbox{ for every chart } (V,x) \mbox{ on } M \} 
\;\;\; .  \]
For all $\varphi \in W^{1,\infty}_\loc(M)$ there is a unique 
$|\nabla \varphi| \in L_{\infty,\loc}(M)$ such that 
\[
|\nabla \varphi| \Big|_V
= \Bigg( \sum_{i,j=1}^d g^{ij} \Big( \frac{\partial}{\partial x^i} \varphi \Big)
                               \Big( \frac{\partial}{\partial x^j} \varphi \Big)
  \Bigg)^{1/2}
\]
for every chart $(V,x)$ on $M$, where 
$\frac{\partial}{\partial x^i} \varphi \in L_{\infty,\loc}(M)$ is defined in the 
natural way.
Then define $W^{1,\infty}(M) = \{ \varphi \in L_\infty(M) : |\nabla \varphi| \in L_\infty(M) \} $.

\begin{prop} \label{pdrum303}
Let $M$ be a connected Riemannian manifold.
If $p,q \in M$ then 
\[
d_M(p\,;q) 
= \sup \{ \psi(p) - \psi(q) : \psi \in W^{1,\infty}_\loc(M), \; 
       |\nabla \psi| \in L_\infty(M) \mbox{ and } \| \, |\nabla \psi| \, \|_\infty \leq 1 \} 
\;\;\; .  \]
\end{prop}
\proof\
`$\leq$'. 
If $q \in M$ define $\psi \colon M \to \Ri$ by $\psi(p) = d_M(p\,;q)$.
Then $\psi \in W^{1,\infty}_\loc(M)$, $\| \, |\nabla \psi| \, \|_\infty \leq 1$ and 
$d_M(p\,;q) \leq \psi(p) - \psi(q)$.

`$\geq$'.
Let $p,q \in M$.
Let $\psi \in W^{1,\infty}_\loc(M)$ with $\| \, |\nabla \psi| \, \|_\infty \leq 1$.
Let $\gamma \colon [0,1] \to M$ be a $C^\infty$-map with $\gamma(0) = p$ and $\gamma(1) = q$.
By regularizing we may assume that $\psi$ is smooth in a neighbourhood of $\gamma([0,1])$.
Then 
\begin{eqnarray*}
|\psi(p) - \psi(q)|
& \leq & \int_0^1 |\dot{\gamma}(t) \, \psi| \, dt  \\
& \leq & \int_0^1 g(\dot{\gamma}(t), \dot{\gamma}(t))^{1/2} \, |\nabla \psi|(\gamma(t)) \, dt
\leq \int_0^1 g(\dot{\gamma}(t), \dot{\gamma}(t))^{1/2} \, dt
\;\;\; . 
\end{eqnarray*}
Minimizing over $\gamma$ gives $|\psi(p) - \psi(q)| \leq d_M(p\,;q)$.\hfill$\Box$

\ruimte

We shall prove that $W^{1,\infty}_\loc(M) = W^{1,\infty}_\loc(N)$ if 
$\capp_M(M \setminus N) = 0$.
For all $s \in [0,\infty)$ we denote by $\ch^s(A)$ the $s$-dimensional Hausdorff measure
of a subset $A$ of $M$ (see \cite{Hei} 8.3).

\begin{prop} \label{pdrum304}
Let $A$ be a subset of a connected Riemannian manifold $M$ of dimension~$d$.
If $\capp(A) = 0$ then $\ch^s(A) = 0$ for all $s \in [0,\infty)$ with $s > d-2$.
\end{prop}
\proof\
For all $n \in \Ni$ there exists a $\varphi_n \in H^1(M)$ such that 
$\varphi_n \geq 1$ on a neighbourhood of $A$ and $\|\varphi_n\|_{H^1(M)} \leq 2^{-n}$.
Set $\varphi = \sum_{n=1}^\infty \varphi_n \in H^1(M)$.
Then for all $m \in \Ni$ it follows that $\varphi \geq m$ on a neighbourhood of $A$.
Hence for all $a \in A$ there exists an $\varepsilon > 0$ such that 
$\langle\varphi\rangle_{a,r} \geq m$ for all $r \in (0,\varepsilon]$, where
$\langle\psi\rangle_{p,r} = |B(p\,;r)|^{-1} \int_{B(p;r)} \psi$ is the 
average of $\psi$ over the ball $B(p\,;r)$
for all $\psi \in L_{1,\loc}(M)$, $p \in M$ and $r > 0$.

Let $a \in A$ and suppose that 
$\limsup_{r \to 0} r^{-s} \int_{B(a;r)} |\nabla \varphi|^2 < \infty$.
Then there exist $r_1 \in (0,1]$ and $M \in \Ri$ such that 
$\int_{B(a;r)} |\nabla \varphi|^2 \leq M \, r^s$ for all $r \in (0,r_1]$.

It follows from \cite{Aub} Theorem~5.14,  that there exists an $r_2 \in (0,r_1]$ such that 
$\exp_a$ is a diffeomorphism from 
$ \{ v \in T_a M : g_a(v,v) < r_2^2 \} $ onto $B(a\,;r_2)$ and 
$d(a\,;\exp_a v) = g_a(v\,;v)^{1/2}$ for all $v \in T_a M$ with $g_a(v,v) < r_2^2$.
Since $T_a M$ is equivalent to $\Ri^d$, it admits a Poincar\'e inequality.
Hence there exists a $c_1 > 0$ such that 
\[
\int_{B(a;r)} |\psi - \langle\psi\rangle_{a,r}|^2 
\leq c_1 \, r^2 \int_{B(a;r)} |\nabla \psi|^2
\]
uniformly for all $\psi \in H^1(B(a\,;r))$ and $r \in (0,r_2]$.
Similarly, there exists a $c_2 > 0$ such that 
$|B(a\,;r)| \geq c_2 \, r^d$ for all $r \in (0,1]$.
Clearly the constants $c_1$ and $c_2$ depend on the point $a$.
Then 
\[
\int_{B(a;r)} |\varphi - \langle\varphi\rangle_{a,r}|^2 
\leq c_1 \, r^2 \int_{B(a;r)} |\nabla \varphi|^2
\leq c_1 \, M \, r^{s+2}
\]
for all $r \in (0,r_2]$.
Therefore
\begin{eqnarray*}
|\langle\varphi\rangle_{a,r} - \langle\varphi\rangle_{a,2r}|
& = & \frac{1}{|B(a\,;r)|} \: \Big| \int_{B(a;r)} \varphi - \langle\varphi\rangle_{a,2r} \Big|  \\
& \leq & \bigg( \frac{1}{|B(a\,;r)|} \, \int_{B(a;r)} |\varphi - \langle\varphi\rangle_{a,2r} |^2 \bigg)^{1/2} \\
& \leq & \bigg( \frac{1}{|B(a\,;r)|} \, \int_{B(a;2r)} |\varphi - \langle\varphi\rangle_{a,2r} |^2 \bigg)^{1/2} \\
& \leq & \Big( 2^{s+2} \, c_1 \, c_2^{-1} M \, r^{s+2-d} \Big)^{1/2}
\end{eqnarray*}
for all $r \in (0,2^{-1} r_2]$.
Since $s+2-d > 0$ one deduces that $(\langle\varphi\rangle_{a,2^{-n}})_n$ is a 
Cauchy sequence.
But for all $m \in \Ni$ there exists an $\varepsilon > 0$ such that 
$\langle\varphi\rangle_{a,r} \geq m$ for all $r \in (0,\varepsilon]$.
This is a contradiction.
Hence $\limsup_{r \to 0} r^{-s} \int_{B(a;r)} |\nabla \varphi|^2 = \infty$.

Let $\delta,\varepsilon > 0$.
There exists a $\gamma > 0$ such that $\int_E |\nabla \varphi|^2 < \varepsilon$ 
for every measurable subset $E$ of $M$ with $|E| < \gamma$.
Since $|A| = 0$ and the measure on $M$ is a Radon measure, there exists 
an open neighbourhood $V$ of $A$ with $|V| < \gamma$.
Set
\[
\cf = \{ B(a\,;r) : a \in A , \; r \in (0,\delta), \; B(a\,;r) \subset V
        \mbox{ and } 
              \int_{B(a;r)} |\nabla \varphi|^2 \geq r^s \} 
\;\;\; .  \]
Then by a basic covering theorem, \cite{Hei} Theorem~1.2, it follows that 
there are $a_1,a_2,\ldots \in A$ and $r_1,r_2,\ldots \in (0,\delta)$ such that 
$B(a_n\,;r_n) \in \cf$ for all $n \in \Ni$, the balls $B(a_n\,;r_n)$ are disjoint
and $\bigcup_{B \in \cf} B \subset \bigcup_{n=1}^\infty B(a_n\,; 5 r_n)$.
Since $A \subset \bigcup_{B \in \cf} B$ it follows that 
\begin{eqnarray*}
\ch^s_{10 \delta}(A)
& \leq & \sum_{n=1}^\infty (5 r_n)^s
\leq 5^s \sum_{n=1}^\infty \int_{B(a_n;r_n)} |\nabla \varphi|^2
\leq 5^s \int_V |\nabla \varphi|^2
< 5^s \varepsilon
\;\;\; .  
\end{eqnarray*}
Hence $\ch^s_{10 \delta}(A) = 0$ and $\ch^s(A) = 0$, as required.\hfill$\Box$

\begin{prop} \label{pdrum305}
Let $N$ be an open subset of a connected Riemannian manifold $M$ of dimension $d$
and suppose that $\ch^{d-1}(M \setminus N) = 0$.
Then $W^{1,\infty}_\loc(N) = W^{1,\infty}_\loc(M)$ and $N$ is connected.
\end{prop}
\proof\
Let $\varphi \in W^{1,\infty}_\loc(N)$.
Using a partition of the unity, normal coordinates and \cite{Hel} Proposition~I.9.10
we may assume that there are $p \in M$ and 
$r > 0$ such that first $\varphi$ is compactly supported in the ball 
$B_M(p\,;r)$, secondly the restriction $\Phi$ of $\exp_p$ to the set 
$X_{2r} =  \{ v \in T_p M : g_p(v,v) < (2r)^2 \} $ is a diffeomorphism of
$X_{2r}$ onto $B_M(p\,;2r)$, thirdly $|v| = d_M(p \,; \exp_p v)$ for all $v \in X_{2r}$
and finally $2^{-1} |v-w| \leq d_M(\exp_p v \,; \exp_p w) \leq 2 |v-w|$
for all $v,w \in X_{2r}$, where $|v| = g_p(v,v)^{1/2}$.
Then 
\[
\ch^{d-1}(X_r \setminus \Phi^{-1}(N \cap B_M(p\,;r)))
= \ch^{d-1}(\Phi^{-1}( (M \setminus N) \cap B_M(p\,;r)))
= 0
\]
where $X_r = \Phi^{-1}(B_M(p\,;r))$.
Moreover, $\varphi \circ \Phi \in W^{1,q}(\Phi^{-1}(N \cap B_M(p\,;r)))$
for all $q \in (1,\infty)$ and $\Phi^{-1}(N \cap B_M(p\,;r))$ is open.
Hence it follows from \cite{AdH} Lemma~9.1.10 that $\varphi \circ \Phi \in W^{1,q}(X_r)$
for all $q \in (1,\infty)$.
Then $\varphi \circ \Phi \in W^{1,\infty}(X_r)$ and 
$\varphi \in W^{1,\infty}(B_M(p\,;r))$.
So $\varphi \in W^{1,\infty}(M)$.

Finally, let $\varphi \colon N \to \{ 0,1 \} $ be a continuous function.
Then $\varphi \in W^{1,\infty}(N) \subset W^{1,\infty}_\loc(M)$.
So $\varphi$ extends to a continuous function on $M$.
But $M$ is connected.
Therefore $\varphi$ is constant and $N$ is connected.\hfill$\Box$

\ruimte

\noindent
{\bf Proof of Theorem~\ref{tdrum302}} \hspace{5pt}\
This easily follows from the last three propositions.\hfill$\Box$

\section{Quasi isometries are isometries} \label{Sdrum5}

In this section we prove that two connected Riemannian manifolds,
which are regular in capacity, are isomorphic if they
have isomorphic open subsets whose complements are polar.
Moreover, we give many useful tools to understand and to work with the $H^1_0$-spaces
defined on Riemannian manifolds.

Let $(M_1,g_1)$ and $(M_2,g_2)$ be Riemannian manifolds. 
We say that
\[ M_1 \simh M_2
\]
if there exist open sets $M_1'\subset M_1$ and $M_2'\subset M_2$  and
an isometry $\tau$ from $M_2'$ onto $M_1'$ such that 
$\capp_{M_1}(M_1\setminus M_1') = 0 = \capp_{M_2}(M_2\setminus M_2')$.

\smallskip

The following theorem is the main theorem in this section.
It shows that the relation $\simh$ defined on connected Riemannian manifolds
determines the manifold.

\begin{thm} \label{tdrum324}
Let $(M_1,g_1)$ and $(M_2,g_2)$ be two connected Riemannian manifolds which are regular
in capacity.
Then
\[
M_1\simh M_2 \Longleftrightarrow (M_1,g_1) \mbox{ and } (M_2,g_2) \mbox{ are isomorphic.}
\]
Explicitly, if $M_1'$ and $M_2'$ are open subsets
of $M_1$ and $M_2$ such that 
$\capp_{M_1}(M_1\setminus M_1') = 0 = \capp_{M_2}(M_2\setminus M_2')$
and $\tau \colon M_2' \to M_1'$ is an isometry,
then there exists an isometry 
$\hat \tau \colon M_2 \to M_1$ such that $\hat \tau|_{M_2'} = \tau$.
\end{thm}

We define the space $H^1_c(M)$ to be 
the set of all $\varphi \in H^1(M)$ such that there exists a compact subset 
$K$ of $M$ with $\varphi = 0$ a.e.\ on $M \setminus K$.

\begin{lemma} \label{ldrum310}
Let $N$ be an open subset of a Riemannian manifold $(M,g)$
such that $|M \setminus N| = 0$.
Then the following are equivalent.
\begin{tabel}
\item \label{ldrum310-1}
$\capp_M(M \setminus N) = 0$.
\item \label{ldrum310-2}
The restriction map $\psi \mapsto \psi|_N$ from $H^1(M)$ into $H^1(N)$
maps $H^1_0(M)$ onto $H^1_0(N)$.
\item \label{ldrum310-3}
$H^1_0(M) 
= \{ \psi|_N : \psi \in H^1_0(M) \} 
= H^1_0(N) 
$.
\end{tabel}
\end{lemma}
\proof\
Clearly \ref{ldrum310-3} is a reformulation of \ref{ldrum310-2}.
`\ref{ldrum310-3}$\Rightarrow$\ref{ldrum310-1}'.
Let $K \subset M$ be a compact set.
There exists a $\psi \in C_c^\infty(M)$ such that 
$\psi|_K \geq 1$.
Then $\psi|_N \in H^1_0(N)$ by assumption.
Let $\varepsilon > 0$.
There exists a $\varphi \in C_c^\infty(N)$ such that 
$\|\psi|_N - \varphi\|_{H^1(N)}^2 < \varepsilon$.
Then $\psi - \varphi \in C_c^\infty(M)$ and 
$\psi - \varphi \geq 1$ on $K \setminus N$.
So 
\[
\capp(K \setminus N) 
\leq \|\psi - \varphi\|_{H^1(M)}^2 
= \|\psi|_N - \varphi\|_{H^1(N)}^2 
< \varepsilon
\;\;\; ,  \]
where we used that $|M \setminus N| = 0$ in the 
equality.
Since $M$ is $\sigma$-compact one deduces that 
$\capp(M \setminus N) = 0$.

`\ref{ldrum310-1}$\Rightarrow$\ref{ldrum310-3}'.
Clearly $\{ \psi|_N : \psi \in H^1_0(M) \} \supset H^1_0(N)$.
Conversely, let $\psi \in C_c^\infty(M)$.
Let $\varepsilon > 0$.
There exists an open neighbourhood of 
$M \setminus N$ and a $\chi \in H^1(M)$ such that 
$0 \leq \chi \leq 1$, $\chi|_U = 1$ and 
$\|\chi\|_{H^1(M)} < \varepsilon$.
Then $(\psi (\one - \chi))|_N \in H^1_c(N) \subset H^1_0(N)$ and 
\[
\|\psi|_N - (\psi (\one - \chi))|_N\|_{H^1(N)}
= \|\psi \, \chi\|_{H^1(M)}
\leq 3 \|\psi\|_{W^{1,\infty}(M)} \, \|\chi\|_{H^1(M)}
\leq 3 \|\psi\|_{W^{1,\infty}(M)} \, \varepsilon
\;\;\; .  \]
So $\psi|_N \in H^1_0(N)$.
Since $C_c^\infty(M)$ is dense in $H^1_0(M)$ and 
$\varphi \mapsto \varphi|_N$ is isometric from 
$H^1(M)$ into $H^1(N)$ the lemma follows.\hfill$\Box$

\ruimte

If $N$ is an open subset of a Riemannian manifold $M$
and $A \subset N$, then $\capp_N(A) \leq \capp_M(A)$.
The next lemma is instrumental to deduce equality 
of the two capacities if $\capp_M(M \setminus N) = 0$.
It is a kind of $L_2$-version of Propositions~\ref{pdrum304} and \ref{pdrum305}.

\begin{lemma} \label{ldrum301}
Let $N$ be an open subset of a manifold $M$.
Suppose that $\capp_M(M \setminus N) = 0$.
Then $H^1(N) = H^1(M)$.
\end{lemma}
\proof\
Let $\varphi \in H^1(N) \cap L_\infty(N)$.
We shall prove that $\varphi \in H^1(M)$.
Let $n \in \Ni$.
Since $\capp_M(M \setminus N) = 0$ there exists a $\psi_n \in H^1(M)$ such that 
$\psi_n \geq 1$ in a neighbourhood of $M \setminus N$ and 
$\|\psi_n\|_{H_1(M)}^2 \leq n^{-1}$.
We may assume that $0 \leq \psi_n \leq 1$.
Then $\varphi - \varphi \, \psi_n \in H^1(N)$.
But $\varphi - \varphi \, \psi_n = \varphi (\one - \psi_n)$ vanishes in 
a neighbourhood of $M \setminus N$.
Therefore we can extend this function by zero
to a function $\varphi_n \in H^1(M)$.
Then 
\[
\|\varphi \, \psi_n\|_{H^1(N)}^2
\leq 2 \|\varphi\|_\infty^2 \|\psi_n\|_{H^1(M)}^2 + 2 \|\varphi\|_{H^1(N)}^2 \|\psi_n\|_\infty^2
\leq  2 \|\varphi\|_\infty^2 + 2 \|\varphi\|_{H^1(N)}^2
\]
for all $n \in \Ni$.
So the sequence $\varphi_1,\varphi_2,\ldots$ is uniformly bounded in $H^1(M)$.
Hence it has a weakly convergent subsequence.
Passing to a subsequence if necessary, there exists a $\tilde \varphi \in H^1(M)$
such that $\lim \varphi_n = \tilde \varphi$ weakly in $H^1(M)$.
Then $\lim \varphi_n = \tilde \varphi$ weakly in $L_2(M)$.
But 
\[
\|\varphi - \varphi_n\|_{L_2(M)}
= \|\varphi \, \psi_n\|_{L_2(M)} 
\leq \|\varphi\|_\infty \, \|\psi_n\|_{L_2(M)} 
\leq \|\varphi\|_\infty \, n^{-1/2}
\]
for all $n \in \Ni$.
So $\lim \varphi_n = \varphi$ strongly and therefore also weakly in $L_2(M)$.
Hence $\tilde \varphi = \varphi$ a.e.\ and $\varphi \in H^1(M)$.

Thus $H^1(N) \cap L_\infty(N) \subset H^1(M)$.
Since $H^1(N) \cap L_\infty(N)$ is dense in $H^1(N)$ 
the lemma follows.\hfill$\Box$

\begin{cor} \label{cdrum320}
If $M$ is a Riemannian manifold and $N \subset M$ is open with 
$\capp_M(M \setminus N)=0$ then $\capp_N(A) = \capp_M(A)$ 
for all $A \subset N$.
\end{cor}

\begin{cor} \label{cdrum321}
Let $M_1$ and $M_2$ be two Riemannian manifolds. 
If $M_1\simh M_2$ and $M_1'$, $M_2'$ and $\tau$ are as in the 
definition of $\simh$ then
$\capp_{M_1}(\tau(A)) = \capp_{M_2}(A)$ for every set $A\subset M_2'$.
\end{cor}
\proof\
One deduces from the previous corollary that 
$\capp_{M_1}(\tau(A)) = \capp_{M_1'}(\tau(A)) 
       = \capp_{M_2'}(A)=\capp_{M_2}(A)=0$.\hfill$\Box$

\ruimte

We emphasize that the next proposition does not require the manifolds
to be regular in capacity.

\begin{prop} \label{pdrum322}
The relation $\simh$ is an equivalence relation. 
\end{prop}
\proof\
The reflexivity and symmetry are trivial.

Let $M_1$, $M_2$ and $M_3$ be three Riemannian manifolds and assume that  
$M_1\simh M_2$ and $M_2\simh M_3$.
Then there exist open $M_1' \subset M_1$, $M_2',M_2'' \subset M_2$ and 
$M_3'' \subset M_3$ and isometries
$\tau \colon M_2' \to M_1$ and $\sigma \colon M_3'' \to M_2''$ such that 
$\capp_{M_1}(M_1 \setminus M_1') = \capp_{M_2}(M_2 \setminus M_2') = 
\capp_{M_2}(M_2 \setminus M_2'') = \capp_{M_3}(M_3 \setminus M_3'') = 0$.
Now let $M_2''' = M_2' \cap M_2''$.
Then $M_2'''$ is open in $M_2$ and
\[ 
\capp_{M_2}(M_2\setminus M_2''')
\leq \capp_{M_2}(M_2\setminus M_2')+\capp_{M_2}(M_2\setminus M_2'')=0
\;\;\; .  \]
Next set $M_1''' = \tau(M_2''') \subset M_1'$
and $M_3''' = \sigma^{-1}(M_2''') \subset M_3'$.
Then $M_1'''$ is open in $M_1$ and $M_3'''$ is open in $M_3$.
Moreover, $\tau|_{M_2'''} \circ \sigma|_{M_3'''}$ is an isometry from
$M_3'''$ onto $M_1'''$.
Since $M_1' \setminus M_1''' = \tau(M_2' \setminus M_2''')$ it follows from
Corollary~\ref{cdrum321} that 
$\capp_{M_1} (M_1' \setminus M_1''') = \capp_{M_2}(M_2' \setminus M_2''') = 0$.
Therefore 
$\capp_{M_1}(M_1 \setminus M_1''') 
   \leq \capp_{M_1}(M_1 \setminus M_1') + \capp_{M_1}(M_1' \setminus M_1''') = 0$.
So $\capp_{M_1}(M_1 \setminus M_1''') = 0$.
It similarly follows that $\capp_{M_3}(M_3 \setminus M_3''') = 0$.
Therefore $M_1 \simh M_3$.\hspace*{10pt}\hfill$\Box$

\ruimte

Also the next proposition does not assume regular in capacity. 
But it overshoots the conclusions in Theorem~\ref{tdrum324}
since the range of $\tilde\tau$ can be bigger than $M_1$.

\begin{prop} \label{pdrum323}
Let $M_1$ and $M_2$ be two connected Riemannian manifolds.
If $M_1\simh M_2$ and if $M_1'$, $M_2'$ and $\tau$ are
as in the definition of $\simh$,
then there exists a distance preserving isomorphism 
$\tilde\tau \colon \widetilde M_2 \to \widetilde M_1$
such that $\tilde\tau|_{M_2'}=\tau$.
\end{prop}
\proof\
The function $\tau|_{M_2'} \colon M_2' \to M_1'$ is distance preserving
with respect to the 
distances $d_{M_2'}$ and $d_{M_1'}$.
Then by Theorem~\ref{tdrum302} the map $\tau|_{M_2'}$ is also an 
distance preserving 
with respect to the induced distances from $M_2$ and $M_1$ on $M_2'$ and $M_1'$.
Since $M_2'$ is dense in $M_2$ and therefore also in $\widetilde M_1$
it follows that there exists a unique distance preserving map 
$\tilde \tau \colon \widetilde M_2 \to \widetilde M_1$
such that $\tilde \tau |_{M_2'} = \tau$.
Similarly there exists a unique distance preserving map
$\tilde \sigma \colon \widetilde M_1 \to \widetilde M_2$
such that $\tilde \sigma |_{M_1'} = \tau^{-1}$.
Then $\tilde \tau \circ \tilde \sigma|_{M_1'}$ is the identity function on $M_1'$.
So by density and continuity $\tilde \tau \circ \tilde \sigma = I_{\widetilde M_1}$
Similarly $\tilde \sigma \circ \tilde \tau = I_{\widetilde M_2}$ and the proposition follows.\hfill$\Box$

\ruimte

Now we are able to prove the main theorem of this section.

\ruimte

\noindent
{\bf Proof of Theorem~\ref{tdrum324} \hspace{5pt}}\
The implication $\Leftarrow$ is trivial.
Suppose that $M_1 \simh M_2$ and let $M_1'$, $M_2'$ and $\tau$ be as in 
the definition of $\simh$.
Let $\tilde \tau$ be as in Proposition~\ref{pdrum323}.
If $\varphi \in H^1_0(M_1) \cap C_0(\widetilde M_1)$ then 
$\varphi \circ \tilde \tau \in C_0(\widetilde M_2)$.
Moreover, $\varphi|_{M_1'} \in H^1_0(M_1')$, so 
$(\varphi \circ \tilde \tau)|_{M_2'} \in H^1_0(M_2')$ and therefore
$(\varphi \circ \tilde \tau)|_{M_2} \in H^1_0(M_2)$ by Lemma~\ref{ldrum310}.
So we can define $V \colon H^1_0(M_1) \cap C_0(\widetilde M_1) \to H^1_0(M_2) \cap C_0(\widetilde M_2)$ by 
$V \varphi = \varphi \circ \tilde \tau$.

Next, let $p \in M_1$.
There exists a $\varphi \in C_c^\infty(M_1)$ such that $\varphi(p) = 1$.
Then $V \varphi \in H^1_0(M_2) \cap C_0(\widetilde M_2)$.
Moreover, $(V \varphi)(\tilde \tau^{-1}(p)) = \varphi(p) = 1$.
Hence $\tilde \tau^{-1}(p) \in M_2$ since $M_2$ is regular in capacity.
Similarly $\tilde \tau(q) \in M_1$ for all $q \in M_2$ since 
$M_1$ is regular in capacity.
Let $\hat \tau = \tilde \tau|_{M_2}$.
Then $\hat \tau$ is a topological 
homeomorphism from $M_2$ onto $M_1$.
It remains to show that $\hat \tau$ and its inverse are smooth and an isometry.

By (\ref{eSdrum1;1}) we can define $U \colon L_2(M_1) \to L_2(M_2)$ by
$U \varphi = \varphi \circ \hat \tau$.
If $\varphi \in H^1_0(M_1)$ then 
$\varphi \circ \tau \in H^1_0(M_2') = H^1_0(M_2)$ by isometry, (\ref{eSdrum1;1})
and Lemma~\ref{ldrum310}.
Therefore $U$ is a bijection from $H^1_0(M_1)$ onto $H^1_0(M_2)$.
Let $h_1$ and $h_2$ be the forms associated to the Dirichlet Laplace--Beltrami operators 
on $M_1$ and $M_2$, with form domains $H^1_0(M_1)$ and $H^1_0(M_2)$.
Then 
\begin{eqnarray*}
h_1(\varphi)
= \int_{M_1} |\nabla \varphi|^2
& = & \int_{M_1'} |\nabla \varphi|^2  \\
& = & \int_{M_2'} |\nabla (\varphi \circ \tau)|^2
= \int_{M_2} |\nabla (\varphi \circ \tau)|^2
= \int_{M_2} |\nabla (\varphi \circ \hat \tau)|^2
= h_2(U \varphi)
\end{eqnarray*}
for all $\varphi \in C_c^\infty(M_1)$.
Since $C_c^\infty(M_1)$ is a core for $h_1$ it follows that 
$h_1(\varphi) = h_2(U \varphi)$ for all $\varphi \in H^1_0(M_1)$.
Hence if $\varphi \in L_2(M_1)$, then $\varphi \in D(\Delta_1)$ if and only if 
$U \varphi \in D(\Delta_2)$, and 
$\Delta_2 U \varphi = U \Delta_1 \varphi$ if both conditions are valid.
Now let $\varphi \in C_c^\infty(M_1)$.
Then $U \varphi \in \bigcap_{n=1}^\infty D(\Delta_2^n) \subset C^\infty(M_2)$
by elliptic regularity.
(See \cite{GT} Theorem~9.11 if $p \neq 1$.
If $p = 1$ first apply a 
Sobolev embedding to embed $L_1$ into a Sobolev space $W^{s,p}$ with $s < 0$ and 
$p > 1$, and then apply \cite{GT} Theorem~9.11.)
So there exists a $\psi \in C^\infty(M_2)$ such that $U \varphi = \psi$ a.e.
But $\varphi = \varphi \circ \hat \tau$ is continuous.
Therefore $\varphi \circ \hat \tau = \psi$ pointwise.
Thus $\varphi \circ \hat \tau$ is smooth for all $\varphi \in C_c^\infty(M_1)$.
Therefore $\hat \tau$ is a $C^\infty$-map from $M_2$ onto $M_1$.
Similarly also $\hat \tau^{-1}$ is a $C^\infty$-map, 
so $\hat \tau$ is a $C^\infty$-diffeomorphism.
Finally, since $\tau$ is an isometry and $M_2'$ is dense in $M_2$ it
follows by continuity that also $\hat \tau$ is an isometry.
This proves Theorem~\ref{tdrum324}.\hfill$\Box$

\section{Lattice homomorphisms} \label{Sdrum2}

In this section we consider lattice homomorphisms between 
$L_p$-spaces on two Riemannian manifolds without the assumption
that the manifolds are regular in capacity.
The aim is to prove that the associated $H^1_0$-spaces are 
equivalent, under the conditions of Theorem~\ref{tdrum101}.

The first step is to use elliptic regularity of the Laplace--Beltrami operator
to reduce to smooth functions.

\begin{lemma} \label{ldrum201}
Let $(M_1,g_1)$ and $(M_2,g_2)$ be two Riemannian manifolds.
Let $p \in [1,\infty)$.
For all $j \in \{ 1,2 \} $ let $\Delta_j$ be the Dirichlet Laplace--Beltrami operator 
on $M_j$ and let $S^{(j)}$ be the associated semigroup on $L_p(M_j)$.
Let  $U \colon L_p(M_1) \to L_p(M_2)$ be a lattice homomorphism such that 
\begin{equation}
U S^{(1)}_t = S^{(2)}_t U
\label{eldrum201;1}
\end{equation}
for all $t > 0$.
Then
\begin{subtabel}
\item \label{ldrum201-1}
$U C_c^\infty(M_1) \subset C^\infty(M_2)$.
\item \label{ldrum201-2}
$U \varphi \geq 0$ for all $\varphi \in C_c^\infty(M_1)$ with $\varphi \geq 0$.
\item \label{ldrum201-3}
$(U \varphi) (U \psi) = 0$ for all $\varphi,\psi \in C_c^\infty(M_1)$ with $\varphi \, \psi = 0$.
\item \label{ldrum201-4}
$\Delta_2 U \varphi = U \Delta_1 \varphi$ for all $\varphi \in C_c^\infty(M_1)$.
\end{subtabel}
\end{lemma}
\proof\
It follows from (\ref{eldrum201;1}) that 
$U D(\Delta_1) \subset D(\Delta_2)$ and $\Delta_2 U \varphi = U \Delta_1 \varphi$
for all $\varphi \in D(\Delta_1)$.
Hence by iteration
$U \bigcap_{n=1}^\infty D(\Delta_1^n) \subset \bigcap_{n=1}^\infty D(\Delta_2^n)$.
But $C_c^\infty(M_j) \subset \bigcap_{n=1}^\infty D(\Delta_j^n) \subset C^\infty(M_j)$ 
for all $j \in \{ 1,2 \} $ by elliptic regularity.
This shows \ref{ldrum201-1} and \ref{ldrum201-4}.
Property~\ref{ldrum201-2} follows since $U$ is a lattice homomorphism.
Moreover, 
$|U \varphi| = U |\varphi|$ for all $\varphi \in L_p(M_1)$.
Hence if $\varphi,\psi \in C_c(M_1)$ and $\varphi \, \psi = 0$ then 
$|\varphi| \wedge |\psi| = 0$ and 
$|U \varphi| \wedge |U \psi| = U |\varphi| \wedge U |\psi|
= U(|\varphi| \wedge |\psi|) = 0$.
Therefore $|(U \varphi) (U \psi)| = |U \varphi| \, |U \psi| = 0$
and $(U \varphi) (U \psi) = 0$.
This implies Property~\ref{ldrum201-3}.\hfill$\Box$



\ruimte

We frequently need the following sufficient condition for point evaluations.

\begin{lemma} \label{ldrum201.3}
Let $M$ be a manifold and $F \colon C_c^\infty(M) \to \Ri$ a positive linear functional
such that 
\begin{equation}
F(\varphi) \, F(\psi) = 0
\mbox{ for all } \varphi,\psi \in C_c^\infty(M) \mbox{ with } \varphi \, \psi = 0
\;\;\; .  
\label{eldrum201.3;1}
\end{equation}
Then there exist $c \in [0,\infty)$ and $p \in M$ such that 
$F(\varphi) = c \, \varphi(p)$ for all $\varphi \in C_c^\infty(M)$.
\end{lemma}
\proof\
Arguing as in \cite{EvG} Corollary~1.8.1 it follows that 
there exists a unique Radon measure $\mu$ on $M$ such that
$F(\varphi) = \int \varphi \, d \mu$ for all $\varphi \in C_c^\infty(M)$.
Then it follows from (\ref{eldrum201.3;1}) that $\mu$ is a point measure.
Hence there exist $p \in M$ and $c \in [0,\infty)$ such that 
$F(\varphi) = c \, \varphi(p)$ for all $\varphi \in C_c^\infty(M)$.\hfill$\Box$

\ruimte

The next proposition is a manifold version of \cite{AB} Theorem~7.22.

\begin{prop} \label{pdrum202}
Let $(M_1,g_1)$ and $(M_2,g_2)$ be two  Riemannian manifolds.
Suppose there exists a 
linear map $U \colon C_c^\infty(M_1) \to C^\infty(M_2)$
such that 
\begin{subtabel}
\item \label{pdrum202-1}
$U \varphi \geq 0$ for all $\varphi \in C_c^\infty(M_1)$ with  $\varphi \geq 0$, and, 
\item \label{pdrum202-2}
$(U \varphi) (U \psi) = 0$ for all $\varphi,\psi \in C_c^\infty(M_1)$ with  $\varphi \, \psi = 0$.
\end{subtabel}
Then there exist an open set $M_2' \subset M_2$, a function
$\tau \colon M_2 \to M_1$ and a function 
$h \colon M_2 \to [0,\infty)$ such that 
$M_2' = \{ q \in M_2 : h(q) > 0 \} $ and 
$U \varphi = h \cdot (\varphi \circ \tau)$ pointwise for all 
$\varphi \in C_c^\infty(M_1)$.
Moreover, the restrictions $\tau|_{M_2'}$ and $h|_{M_2'}$ are both $C^\infty$.
\end{prop}
\proof\
Set $M_2' = \{ q \in M_2 : (U \varphi)(q) \neq 0 \mbox{ for some } \varphi \in C_c^\infty(M_1) \} $.
Then $M_2'$ is open.
Let $q \in M_2$.
Then the map $\varphi \mapsto (U \varphi)(q)$ from $C_c^\infty(M_1)$ into 
$\Ri$ is linear, positive and satisfies (\ref{eldrum201.3;1}).
Hence it follows from Lemma~\ref{ldrum201.3} that 
there exist $\tau(q) \in M_1$ and $h(q) \in [0,\infty)$ such that 
\[
(U \varphi)(q)
= h(q) \, \varphi(\tau(q))
\]
for all $\varphi \in C_c^\infty(M_1)$.
So one obtains functions $\tau \colon M_2 \to M_1$ and 
$h \colon M_2 \to [0,\infty)$.
Moreover, $M_2' = \{ q \in M_2 : h(q) > 0 \} $.
It remains to show the smoothness of the restrictions of $\tau$ and $h$ 
to the set $M_2'$.

First we show that the function $\tau|_{M_2'}$ is continuous.
Otherwise there are $q,q_1,q_2,\ldots \in M_2'$
and $\varepsilon > 0$ 
such that $\lim q_n = q$ and $d(\tau(q_n),\tau(q)) \geq \varepsilon$ for all $n \in \Ni$,
where $d$ is a distance on $M_1$.
There exists a $\varphi \in C_c^\infty(M_1)$ such that  $\varphi(\tau(q)) = 1$ and 
$\varphi(p) = 0$ for all $p \in M_1$ with $d(p,\tau(q)) > \varepsilon$.
Then $U \varphi$ is continuous, so 
$h(q) = (U \varphi)(q) = \lim (U \varphi)(q_n) = 0$, which is a contradiction.

Secondly, let $\Omega$ be a relatively compact open subset of $M_2'$ with $\overline{\Omega} \subset M_2'$.
Then $\tau(\overline \Omega)$ is compact, so there is a 
$\psi \in C_c^\infty(M_1)$ such that $\psi|_{\tau(\Omega)} = 1$.
Then $h|_\Omega = (U \psi)|_\Omega$ is a $C^\infty$-function.
Hence $h|_{M_2'}$ is a $C^\infty$-function from $M_2'$ into $(0,\infty)$.
Then $(\varphi \circ \tau)|_\Omega = (h^{-1} \cdot U \varphi)|_\Omega \in C^\infty(\Omega)$ for all 
$\varphi \in C_c^\infty(M_1)$ and $\tau|_{M_2'}$ is a $C^\infty$-function.\hfill$\Box$

\ruimte

Using the fact that $U$ intertwines with the Laplace--Beltrami operators
implies that $\tau$ is almost an isometry.

\begin{prop} \label{pdrum202.3}
Let $(M_1,g_1)$ and $(M_2',g_2)$ be two  Riemannian manifolds with 
Dirichlet Laplace--Beltrami operators $\Delta_1$ and $\Delta_2$.
Let $\tau \colon M_2' \to M_1$ be a $C^\infty$-map and 
$h \colon M_2' \to (0,\infty)$ a $C^\infty$-function.
Define the map $U \colon C_c^\infty(M_1) \to C^\infty(M_2')$ by
$U \varphi = h \cdot (\varphi \circ \tau)$.
Suppose that 
$\Delta_2 U = U \Delta_1$.
Then 
\begin{equation}
g_1|_{\tau(q)}(\alpha,\beta) = g_2|_q(\tau^*(\alpha),\tau^*(\beta))
\label{epdrum202.3;4}
\end{equation}
for all $q \in M_2'$
and $\alpha,\beta \in T^*_{\tau(q)} M_1$.
In particular, $\dim M_2' \geq \dim M_1$.

If, in addition, $\dim M_1 = \dim M_2'$ then $\tau$ is locally an isometry and 
$h$ is locally constant.
\end{prop}
\proof\
It follows from the identity $\Delta_2 U = U \Delta_1$ that 
\begin{equation}
\Delta_2( h \cdot (\varphi \circ \tau))
= h \cdot ( (\Delta_1 \varphi) \circ \tau)
\label{epdrum202;2}
\end{equation}
on $M_2'$ for all $\varphi \in C_c^\infty(M_1)$.
Let $q \in M_2'$.
There exists a chart $(V,x)$ on $M_1$ such that 
$\tau(q) \in V$ and $x(\tau(q)) = 0$.
Let $\Omega \subset M_1$ be a relatively compact subset such that 
$\tau(q) \in \Omega \subset \overline\Omega \subset V$.
Let $d_1 = \dim M_1$ and $d_2 = \dim M_2'$.
Let $\lambda_1,\ldots,\lambda_{d_1} \in \Ri$.
For all $t > 0$ there exists a $\varphi_t  \in C_c^\infty(M_1)$
such that 
\[
\varphi_t |_\Omega 
= e^{t \sum_{k=1}^{d_1} \lambda_k x^k} |_\Omega
\;\;\, .  \]
Since
\[
\Delta_1
= - \sum_{i,j=1}^{d_1} \frac{1}{\sqrt{g_1}} \, \frac{\partial}{\partial x^i} \, 
           g_1^{ij} \, \sqrt{g_1} \, \frac{\partial}{\partial x^j}
\]
on $V$ it follows that 
\[
\Delta_1 \varphi_t
= - \sum_{i,j=1}^{d_1} t^2 g_1^{ij} \, \lambda_i \, \lambda_j \, \varphi_t
   + t \, \frac{\lambda_j}{\sqrt{g_1}} \, \varphi_t \, 
     \frac{\partial}{\partial x^i} ( g_1^{ij} \, \sqrt{g_1} )
\]
on $\Omega$.
Hence 
\[
\lim_{t \to \infty} t^{-2} \Big( h \cdot ( (\Delta_1 \varphi_t) \circ \tau) \Big) (q)
= - h(q) \sum_{i,j=1}^{d_1} g_1^{ij}|_{\tau(q)} \, \lambda_i \, \lambda_j
\;\;\, .  \]
Next, let $(W,y)$ be a chart on $M_2'$ such that $q \in W$.
Then it follows similarly that 
\begin{eqnarray*}
\lim_{t \to \infty} t^{-2} \Big( \Delta_2( h \cdot (\varphi_t \circ \tau)) \Big) (q)
& = & - \sum_{k,l=1}^{d_2} h(q) \, g_2^{kl}|_q
   \Big( \frac{\partial}{\partial y^k} \sum_{i=1}^{d_1} \lambda_i x^i \circ \tau \Big)(q)
   \Big( \frac{\partial}{\partial y^l} \sum_{j=1}^{d_1} \lambda_j x^j \circ \tau \Big)(q)  \\
& = & - \sum_{i,j=1}^{d_1} h(q) \, g_2|_q(\tau^*(dx^i),\tau^*(dx^j)) \, \lambda_i \, \lambda_j
\;\;\, . 
\end{eqnarray*}
But then (\ref{epdrum202;2}) gives
\[
\sum_{i,j=1}^{d_1} g_1^{ij}|_{\tau(q)} \, \lambda_i \, \lambda_j
= \sum_{i,j=1}^{d_1} g_2|_q(\tau^*(dx^i),\tau^*(dx^j)) \, \lambda_i \, \lambda_j
\]
for all $\lambda_1,\ldots,\lambda_{d_1} \in \Ri$ and 
$g_1^{ij}|_{\tau(q)} = g_2|_q(\tau^*(dx^i),\tau^*(dx^j))$ for all 
$i,j \in \{ 1,\ldots,d_1 \} $.
Hence 
$g_1|_{\tau(q)}(\alpha,\beta) = g_2|_q(\tau^*(\alpha),\tau^*(\beta))$ for all
and $\alpha,\beta \in T^*_{\tau(q)} M_1$.
In particular, $\tau^*$ is injective and $d_2 \geq d_1$.

Finally suppose that $d_1 = d_2$.
Since $\tau^*|_{\tau(q)}$ is injective for all $q \in M_2'$ it follows that 
$\tau^*|_{\tau(q)}$ is bijective for all $q \in M_2'$.
Hence $\tau$ is locally a $C^\infty$-diffeomorphism.
Moreover, $\tau^*|_{\tau(q)}$ is a orthogonal map and therefore also 
$\tau_*|_q$ is an orthogonal map for all $q \in M_2'$.
In particular, $\tau$ is locally an isometry.

It then also follows that 
$\Delta_2(\varphi \circ \tau)
= (\Delta_1 \varphi) \circ \tau$
on $M_2'$ for all $\varphi \in C_c^\infty(M_1)$.
If $q \in M_2'$ and $\Omega$ is a relatively compact open subset of $M_2'$ 
with $q \in \Omega \subset \overline \Omega \subset M_2'$, and if one chooses 
$\varphi \in C_c^\infty(M_1)$ such that $\varphi|_{\tau(\Omega)} = 1$,
then it follows from (\ref{epdrum202;2}) that 
$(\Delta_2 h)(q) = 0$.
Then for all $\varphi \in C_c^\infty(M_1)$ one deduces from (\ref{epdrum202;2}) that
\begin{eqnarray*}
h \cdot ( (\Delta_1 \varphi) \circ \tau)
& = & \Delta_2( h \cdot (\varphi \circ \tau))  \\
& = & (\Delta_2 h) \cdot (\varphi \circ \tau) + 2 (\nabla_2 h) \cdot (\nabla_2 (\varphi \circ \tau)) 
          + h \cdot \Delta_2(\varphi \circ \tau)  \\
& = & 2 (\nabla_2 h) \cdot (\nabla_2 (\varphi \circ \tau)) + h \cdot ( (\Delta_1 \varphi) \circ \tau)
\end{eqnarray*}
on $M_2'$.
So $(\nabla_2 h) \cdot (\nabla_2 (\varphi \circ \tau)) = 0$ on $M_2'$.
Since $\tau$ is locally a diffeomorphism it follows that $\nabla_2 h = 0$ and 
$h$ is locally constant.
This completes the proof of Proposition~\ref{pdrum202.3}.\hfill$\Box$

\ruimte

In the next lemmas we consider injectivity and the density of the range of $U$.

\begin{lemma} \label{ldrum202.5}
Assume the hypothesis of Lemma~$\ref{ldrum201}$.
Moreover, assume  $U \neq 0$, the manifold $M_1$ is connected 
and  the manifolds
have equal dimension.

Then there exist open sets $M_1' \subset M_1$ and $M_2' \subset M_2$, a  
map $\tau \colon M_2 \to M_1$ and a bounded function $h \colon M_2 \to [0,\infty)$
such that $M_2' = \{ q \in M_2 : h(q) > 0 \} $, $M_1' = \tau(M_2')$,
$\tau|_{M_2'}$ is a local isometry, $h|_{M_2'}$ is a $C^\infty$-function and 
\[
U \varphi = h \cdot (\varphi \circ \tau)
\]
pointwise for all $\varphi \in C_c^\infty(M_1)$ and a.e.\ for all $\varphi \in L_p(M_1)$.
Moreover, $U$ is injective and $|M_1 \setminus M_1'| = 0$.
\end{lemma}
\proof\
Since $C_c^\infty(M_1)$ is dense in $L_p(M_1)$ and 
$U \neq 0$ it follows that the restriction of $U$
to $C_c^\infty(M_1)$ does not vanish.
So we can apply Lemma~\ref{ldrum201} and Proposition~\ref{pdrum202}.
We use the notation of Proposition~\ref{pdrum202}.
Set $M_1' = \tau(M_2')$.
It follows from the inverse function theorem that $M_1'$ is 
open in $M_1$.
Moreover, since $M_2' \neq \emptyset$ it follows that 
$M_1' \neq \emptyset$ and $|M_1'| \neq 0$.

Now let $\varphi \in L_p(M_1)$.
Since $C_c^\infty(M_1)$ is dense in $L_p(M_1)$
there exists a sequence $\varphi_1,\varphi_2,\ldots \in C_c^\infty(M_1)$
such that $\lim \varphi_n = \varphi$ in $L_p(M_1)$.
Since $U$ is continuous it follows that $\lim U \varphi_n = U \varphi$
in $L_p(M_2)$.
Passing to subsequences, if necessary, we may assume that 
$\lim \varphi_n = \varphi$ a.e.\ and 
$\lim U \varphi_n = U \varphi$ a.e.
But since $M_1'$ is $\sigma$-compact and $\tau|_{M_2'}$ is locally an 
isometry it follows that $\tau^{-1}(N)$ is a null-set in $M_2'$ 
for every null-set $N$ in $M_1'$.
Therefore $\lim \varphi_n \circ \tau = \varphi \circ \tau$ a.e.\ on $M_2'$.
Then 
\[
U \varphi 
= \lim U \varphi_n
= \lim h \cdot (\varphi_n \circ \tau)
= h \cdot (\varphi \circ \tau)
\]
a.e.\ on $M_2'$.
In addition one has $U \varphi = \lim U \varphi_n = 0$ a.e.\ on $M_2 \setminus M_2'$.
So $U\varphi=h\cdot (\varphi \circ\tau)$ a.e.\ on $M_2$.

Next we show that $U$ is injective.
Let $\varphi \in L_2(M_1)$ and suppose that $U \varphi = 0$.
Then $U |\varphi| = |U \varphi| = 0$, so we may assume that $\varphi \geq 0$.
Fix $t > 0$.
Then 
\[
0 
= S^{(2)}_t U \varphi
= U S^{(1)}_t \varphi
= h \cdot (S^{(1)}_t \varphi) \circ \tau
\]
a.e.\ on $M_2'$.
Since also $\tau$ maps $M_2'$-null-sets into $M_1'$-null-sets 
and $h > 0$ pointwise, one deduces
that $S^{(1)}_t \varphi = 0$ a.e.\ on $M_1'$.
But $M_1$ is connected and if $\varphi \neq 0$ then $(S^{(1)}_t \varphi)(p) > 0$
for all $p \in M_1$ and in particular for all $p \in M_1'$.
So $\varphi = 0$ and $U$ is injective.
It is obvious that this implies that $|M_1 \setminus M_1'| = 0$.

Finally we show that $h$ is bounded by $\|U\|$.
Let $q \in M_2'$ and $\varepsilon > 0$.
Since $\tau|_{M_2'}$ is locally an isometry and $h|_{M_2'}$ is continuous there exists an open 
neighbourhood $\Omega$ of $q$ in $M_2'$ such that $\tau|_\Omega \colon \Omega \to \tau(\Omega)$
is an isometry and $h|_\Omega \geq (1 - \varepsilon) h(q)$.
Fix $\varphi \in C_c^\infty(\tau(\Omega))$ with $\varphi \neq 0$.
Then 
\begin{eqnarray*}
(1 - \varepsilon) \, h(q) \|\varphi \circ \tau\|_{L_p(\Omega)}
& \leq & \|h \cdot (\varphi \circ \tau)\|_{L_p(\Omega)}
\leq \|U \varphi\|_{L_p(M_2)}  \\
& \leq & \|U\| \, \|\varphi\|_{L_p(M_1)}
\leq \|U\| \, \|\varphi\|_{L_p(\tau(\Omega))}
= \|U\| \, \|\varphi \circ \tau\|_{L_p(\Omega)}
\end{eqnarray*}
where we used (\ref{eSdrum1;1}) in the last step.
So $h(q) \leq \|U\|$.\hfill$\Box$

\begin{lemma} \label{ldrum202.6}
Assume the hypothesis of Lemma~$\ref{ldrum201}$.
Let $M_2'$ be an open subset of $M_2$,
let $\tau \colon M_2 \to M_1$ and 
$h \colon M_2 \to [0,\infty)$ be two functions such that 
$M_2' = \{ q \in M_2 : h(q)> 0\}$ and 
the restrictions $\tau|_{M_2'}$ and $h|_{M_2'}$ are both $C^\infty$-maps.
Suppose that $U \varphi = h \cdot (\varphi \circ \tau)$ for all 
$\varphi \in C_c^\infty(M_1)$.

Then the following are equivalent.
\begin{tabel}
\item \label{ldrum202.6-1}
$U L_p(M_1)$ is dense in $L_p(M_2)$.
\item \label{ldrum202.6-2}
$U C_c^\infty(M_1)$ is dense in $L_p(M_2)$.
\item \label{ldrum202.6-3}
For every pair of disjoint measurable
subsets $A_1$ and $A_2$ in $M_2$  with $0 < |A_1|, |A_2| < \infty$
the functionals
\[
\varphi \mapsto \int_{A_1} U \varphi
\;\;\;\; \mbox{and} \;\;\;\;
\varphi \mapsto \int_{A_2} U \varphi
\]
from $C_c^\infty(M_1)$ into $\Ri$ are linearly independent.
\item \label{ldrum202.6-4}
The map $\tau|_{M_2'}$ is injective  and $|M_2 \setminus M_2'| = 0$.
\end{tabel}
Moreover, these conditions imply that the dimensions of $M_1$ and $M_2$ are equal.
\end{lemma}
\proof\
Clearly \ref{ldrum202.6-1}$\Leftrightarrow$\ref{ldrum202.6-2}$\Rightarrow$\ref{ldrum202.6-3}.

Next we show that \ref{ldrum202.6-3} or \ref{ldrum202.6-4} 
implies that $\dim M_1 = \dim M_2$.
Obviously $U \neq 0$, so $M_2' \neq 0$ and $d_2 \geq d_1$ by Proposition~\ref{pdrum202.3},
where $d_1 = \dim M_1$ and $d_2 = \dim M_2$.
Fix $q \in M_2'$ and set $p = \tau(q)$.
Let $(V,x)$ be a chart on $M_1$ with $p \in V$ and $(W,y)$ a chart on $M_2'$ 
with $q \in W$.
Let $V' = x(V)$, $W' = y(W)$, $p' = x(p)$ and $q' = y(q)$.
Define the $C^\infty$-map $f \colon W' \to V'$ by 
$f = x \circ \tau \circ y^{-1}$.
Then it follows from (\ref{epdrum202.3;4}) that 
$(Df)(q')$ has maximal rank.
Suppose that $k = d_2 - d_1 > 0$.
Then it follows from the inverse function theorem that 
there exist open $W'' \subset \Ri^{d_2}$, $\delta > 0$ and 
a $C^\infty$-diffeomorphism $F \colon W'' \to B(p'\,;3 \delta) \times (-3\delta,3\delta)^k$
such that $q \in W'' \subset W'$, 
$B(p'\,;3\delta) \subset V'$, $F(q') = (p',0)$ and $f(G(u,v)) = u$
for all $(u,v) \in B(p'\,;3 \delta) \times (-3\delta,3\delta)^k$, where $G = F^{-1}$.

In particular, $f$ is not injective. 
This contradicts the injectivity of $\tau$ in \ref{ldrum202.6-4}.

In order to obtain a contradiction with Condition~\ref{ldrum202.6-3}
we proceed as follows.
Define the $C^\infty$-function $\Phi \colon B(p'\,;3 \delta) \times (-3\delta,3\delta)^k \to (0,\infty)$
by $\Phi(u,v) = ((h \, \sqrt{g_2}) \circ y^{-1} \circ G)(u,v) \, |(JG)(u,v)|$,
where $JG$ denotes the Jacobian determinant of $G$.
If $A \subset y^{-1}(W'')$ is measurable then 
\begin{equation}
\int_A U \varphi
= \int_{(F \circ y)(A)} \Phi(u,v) \, (\varphi \circ x^{-1})(u) \, d(u,v)
\label{eldrum202.6;4}
\end{equation}
for all $\varphi \in C_c^(M_1)$.
By compactness, there are $m,M > 0$ such that 
$m \leq \Phi(u,v) \leq M$ for all $(u,v) \in \overline{B(p'\,;2\delta)} \times [-2\delta,2\delta]^k$.
Let $u \in \overline{B(p'\,;2\delta)}$.
Then 
\[
\int_{[-2\delta,2\delta]^{k-1} \times [-\delta \, m \, M^{-1},0)} \Phi(u,v) \, dv
\leq (4\delta)^{k-1} \, \delta \, m
\leq \int_{[-2\delta,2\delta]^{k-1} \times (0,\delta]} \Phi(u,v) \, dv
\;\;\; .  \]
So there exists a $\beta(u) \in (0,\delta]$ such that
\[
\int_{[-2\delta,2\delta]^{k-1} \times [-\delta \, m \, M^{-1},0)} \Phi(u,v) \, dv
= \int_{[-2\delta,2\delta]^{k-1} \times (0,\beta(u)]} \Phi(u,v) \, dv
\;\;\; .  \]
Now choose
\[
A_1 = (y^{-1} \circ F^{-1}) \Big( \overline{B(p',;2\delta)} \times [-2\delta,2\delta]^{k-1} \times [-\delta \, m \, M^{-1},0) \Big)
\]
and 
\[
A_2 = \{ (y^{-1} \circ F^{-1})(u,v) : u \in \overline{B(p'\,;2\delta)} \mbox{ and } 
                 v \in [-2\delta,2\delta]^{k-1} \times (0,\beta(u)] \} 
\;\;\; .  \]
Then $A_1 \cap A_2 = \emptyset$, $0 < |A_1|,|A_2| < \infty$ and 
$\int_{A_1} U \varphi = \int_{A_1} U \varphi$ for all $\varphi \in C_c^\infty(M_1)$ by 
(\ref{eldrum202.6;4}) and the choice of $\beta(u)$.
This contradicts \ref{ldrum202.6-3}.
So $d_1 = d_2$.
Thus all four conditions imply that the dimensions of $M_1$ and $M_2$ are equal.

`\ref{ldrum202.6-3}$\Rightarrow$\ref{ldrum202.6-4}'.
Clearly \ref{ldrum202.6-3} implies that $|M_2 \setminus M_2'| = 0$.
Suppose \ref{ldrum202.6-3} and $\tau$ is not injective.
Then there are $q_1,q_2 \in M_2'$ such that $\tau(q_1) = \tau(q_2)$ and $q_1 \neq q_2$.
Since $\dim M_1 = \dim M_2$ it follows from Proposition~\ref{pdrum202.3} that 
$\tau|_{M_2'}$ is locally an isomorphism.
There are open connected relative compact $\Omega_1, \Omega_2 \subset M_2'$ such that 
$\Omega_1 \cap \Omega_2 = \emptyset$ and for all $j \in \{ 1,2 \} $ one has 
$q_j \in \Omega_j$ and $\tau_j = \tau|_{\Omega_j} \colon \Omega_j \to \tau(\Omega_j)$ is an
isometry.
Since $\Omega_j$ is connected there is a $c_j \in (0,\infty)$ such that 
$h|_{\Omega_j} = c_j$.
Without loss of generality we may assume that $\tau_1(\Omega_1) = \tau_2(\Omega_2)$.
Then for all $\varphi \in C_c^\infty(M_1)$ one has 
$(U \varphi)(q) = h(q) \, \varphi(\tau_1(q)) = c_1 \, \varphi(\tau_1(q))$
for all $q \in \Omega_1$.
So
\[
\int_{\Omega_1} U \varphi
= c_1 \int_{\Omega_1} \varphi \circ \tau_1
= c_1 \int_{\tau_1(\Omega_1)} \varphi
\]
where we used (\ref{eSdrum1;1}).
Similarly
\[
\int_{\Omega_1} U \varphi
= c_2 \int_{\tau_2(\Omega_2)} \varphi
\]
for all $\varphi \in C_c^\infty(M_1)$.
But $\tau_1(\Omega_1) = \tau_2(\Omega_2)$.
This contradicts the independence of the functionals.

`\ref{ldrum202.6-4}$\Rightarrow$\ref{ldrum202.6-2}'.
Since $|M_2 \setminus M_2'| = 0$ the space $C_c^\infty(M_2')$ is dense in $L_p(M_2)$.
Therefore it suffices to show that $C_c^\infty(M_2') \subset U C_c^\infty(M_1')$.

Using again that $\dim M_1 = \dim M_2$ it follows from Proposition~\ref{pdrum202.3} that 
$\tau|_{M_2'}$ is locally an isomorphism and $h_{M_2'}$ is locally constant.
If $\psi \in C_c^\infty(M_2')$ and there exists an open connected set $\Omega$ in 
$M_2'$ such that $\supp \psi \subset \Omega$, then $h$ is constant on $\Omega$, say $c$,
and $\varphi = c^{-1} \psi \circ (\tau|_{M_2'})^{-1} \in C_c^\infty(M_1') \subset C_c^\infty(M_1)$
satisfies $U \varphi = \psi$.
Then the general case follows by a partition of the unity.\hfill$\Box$

\ruimte

For open subsets in $\Ri^d$ the surjectivity of $U$ follows from the fact that $U \neq 0$
(see \cite{Are4}, Theorem~2.1). 
In general the condition $U \neq 0$ is not sufficient to 
establish the surjectivity of $U$.

\begin{exam} \label{xdrum501}
Let $S^1 = \{ z \in \Ci : |z| = 1 \} $.
Let $g_1$ be the Riemannian metric on $M_1 = S^1$ such that 
$(g_1)|_{e^{i \theta}}(\frac{\partial}{\partial x^1} |_{e^{i \theta}}, \frac{\partial}{\partial x^1}|_{e^{i \theta}}) = 1$
for each $\theta \in \Ri$, where $(V,x)$ is a chart on $S^1$ such that $V$ is an open neighbourhood
of $e^{i \theta}$, $\theta \in x(V)$ and $x^{-1}(\xi) = e^{i \xi}$ for all $\xi \in x(V)$.
Set $M_2 = S_1$ and choose the Riemannian metric $g_2$ on $M_2$ 
by $g_2 = 4 g_1$.

Define $U \colon L_2(M_1) \to L_2(M_2)$ by 
\[
(U \varphi)(z) = \varphi(z^2)
\;\;\; .  \]
Then $U$ is a lattice homomorphism, $U \neq 0$ and 
$U S^{(1)}_t = S^{(2)}_t U$ for all $t > 0$, where
$S^{(j)}$ is the semigroup 
on $L_2(M_j)$ generated by the Dirichlet Laplace--Beltrami operator on $M_j$ 
for all $j \in \{ 1,2 \} $.
Moreover, $M_1$ and $M_2$ are regular in capacity.
But the Riemannian manifolds $(M_1,g_1)$ and $(M_2,g_2)$ are 
not isomorphic.
\end{exam}

We combine the previous results.

\begin{prop} \label{pdrum203}
Let $(M_1,g_1)$ and $(M_2,g_2)$ be two connected Riemannian manifolds.
Let $p \in [1,\infty)$.
For all $j \in \{ 1,2 \} $ let $\Delta_j$ be the Dirichlet Laplace--Beltrami operator 
on $M_j$ and let $S^{(j)}$ be the associated semigroup on $L_p(M_j)$.
Let  $U \colon L_p(M_1) \to L_p(M_2)$ be a lattice homomorphism such that 
$U L_p(M_1)$ is dense in $L_p(M_2)$ and 
\begin{equation}
U S^{(1)}_t = S^{(2)}_t U
\label{epdrum203;1}
\end{equation}
for all $t > 0$.

Then $U$ is an order isomorphism
and there exist open connected sets $M_1' \subset M_1$ and $M_2' \subset M_2$, a  
map $\tau \colon M_2 \to M_1$ and a constant $c > 0$ such that $M_1' = \tau(M_2')$,
$\tau|_{M_2'} \colon M_2' \to M_1'$ is an isometry, and 
\[
U \varphi = c \, \one_{M_2'} \cdot (\varphi \circ \tau)
\]
pointwise for all $\varphi \in C_c^\infty(M_1)$ and a.e.\ for all $\varphi \in L_p(M_1)$.
Moreover, 
$\capp(M_1 \setminus M_1') = \capp(M_2 \setminus M_2') = 0$
and for all $\tilde p \in [1,\infty)$ there exists an order isomorphism $\widetilde U$ such that 
such that 
$\widetilde U \varphi = U \varphi$ for all $\varphi \in L_p(M_1) \cap L_{\tilde p}(M_1)$.
Finally, if $p = 2$ then $U$ maps $H^1_0(M_1)$ continuously into $H^1_0(M_2)$.
\end{prop}
\proof\
We use the notation as above.
Let $\tilde p \in [1,\infty)$.
Then the map 
$\widetilde U \colon L_{\tilde p}(M_1) \to L_{\tilde p}(M_2)$
defined by $\widetilde U \varphi = h \cdot (\varphi \circ \tau)$ (a.e.) for all 
$\varphi \in L_{\tilde p}(M_1)$ is well defined
since $h$ is bounded and $\tau^{-1}(N)$ is a null-set in $M_2'$ for every null-set $N$
in $M_1'$.
It is a lattice homomorphism and is consistent with $U$.
Moreover, $U L_{\tilde p}(M_1)$ is dense in $L_{\tilde p}(M_2)$ by 
Proposition~\ref{ldrum202.6} \ref{ldrum202.6-4}$\Rightarrow$\ref{ldrum202.6-1}.

Therefore, for the remainder of the proof we may assume that $p = 2$.
Then one deduces from (\ref{epdrum203;1}) that 
$(I + \Delta_2)^{-1/2} U = U (I + \Delta_1)^{-1/2}$.
Hence 
\begin{eqnarray*}
U H^1_0(M_1)
& = & U (I + \Delta_1)^{-1/2} L_2(M_1)  \\
& = & (I + \Delta_2)^{-1/2} U L_2(M_1)
\subset (I + \Delta_2)^{-1/2} L_2(M_2)
= H^1_0(M_2)
\;\;\; .  
\end{eqnarray*}
Then by the closed graph theorem the restriction of $U$
to $H^1_0(M_1)$ is a continuous map from $H^1_0(M_1)$ into $H^1_0(M_2)$.
Next $U L_2(M_1)$ is dense in $L_2(M_2)$ and 
$(I + \Delta_2)^{-1/2}$ is continuous from $L_2(M_2)$ onto $H^1_0(M_2)$.
So $U H^1_0(M_1)$ is dense in $H^1_0(M_2)$.
Therefore $U C_c^\infty(M_1)$ is dense in $H^1_0(M_2)$.

Now suppose $\capp(M_2 \setminus M_2') > 0$.
There exist compact subsets $K_1 \subset K_2 \subset \ldots$ of $M_2$ 
such that $M_2 = \bigcup_{n=1}^\infty K_n$.
Then $\capp(M_2 \setminus M_2') = \lim_{n \to \infty} \capp(K_n \setminus M_2')$
by \cite{BH} Proposition~8.1.3c.
Hence there exists an $n \in \Ni$ such that $\capp(K_n \setminus M_2') > 0$.
Let $\psi \in C_c^\infty(M_2)$ be such that 
$\psi|_{K_n} = 1$.
Since $U C_c^\infty(M_1)$ is dense in $H^1_0(M_2)$ there exists a 
$\varphi \in C_c^\infty(M_1)$ such that 
$\|\psi - U \varphi\|_{H^1(M_2)}^2 < \capp(K_n \setminus M_2')$.
Then $\psi - U \varphi = 1$ on $K_n \setminus M_2'$ by definition of $M_2'$.
Moreover, $\psi - U \varphi$ is continuous and 
$\psi - U \varphi \in H^1_0(M_2) \subset H^1(M_2)$.
Hence $\capp(K_n \setminus M_2') \leq \|\psi - U \varphi\|_{H^1(M_2)}^2$.
This is a contradiction.
Hence $\capp(M_2 \setminus M_2') = 0$.

This allows to apply Theorem~\ref{tdrum302} to deduce that $M_2'$ is connected.
Then $h|_{M_2'}$ is constant, say $c > 0$.
It follows from (\ref{eSdrum1;1}) that (for any $p \in [1,\infty)$) the map 
$U$ is an isometry between $L_p$-spaces and since the range is dense, it is surjective.

Finally, $H^1_0(M_2') = \{ \psi|_{M_2'} : \psi \in H^1_0(M_2) \} $ 
by Lemma~\ref{ldrum310}, since $\capp(M_2 \setminus M_2') = 0$.
But $H^1_0(M_2) = U H^1_0(M_1)$.
In addition, $U H^1_0(M_1') = H^1_0(M_2')$ since 
$\tau|_{M_2'} \colon M_2' \to M_1'$ is an isometry.
Therefore 
\[
U H^1_0(M_1')
= \{ (U \varphi)|_{M_2'} : \varphi \in H^1_0(M_1) \}
= \{ U (\varphi|_{M_1'}) : \varphi \in H^1_0(M_1) \}
\;\;\; .  \]
Hence $H^1_0(M_1') = \{ \varphi|_{M_1'} : \varphi \in H^1_0(M_1) \} $.
Using Lemma~\ref{ldrum310} again one deduces that 
$\capp(M_1 \setminus M_1') = 0$.\hfill$\Box$

\ruimte

Now we are able to prove the main theorem of this paper.

\ruimte

\noindent
{\bf Proof of Theorem~\ref{tdrum101} \hspace{5pt}}\
The implication \ref{tdrum101-1}$\Rightarrow$\ref{tdrum101-3} is trivial.
But if Condition~\ref{tdrum101-3} is valid then Proposition~\ref{pdrum203}
implies that $M_1 \simh M_2$. 
Hence the Riemannian manifolds $(M_1,g_1)$ and $(M_2,g_2)$ are isomorphic by Theorem~\ref{tdrum324}.
Moreover, there exist $c > 0$ and a isometry $\tau \colon M_2 \to M_1$ such that 
$U \varphi = c \, \varphi \circ \tau$ for all $\varphi \in L_p(M_1)$.
\hfill$\Box$

\ruimte

In fact, it follows from the proof that under 
Condition~\ref{tdrum101-3} it follows that $M_2' = M_2$ and $M_1' = M_1$.
Therefore $U \varphi = c \, \varphi \circ \tau$ for all $\varphi \in L_p(M_1)$.

\section{Regularity in capacity} \label{Sdrum4}

The purpose of this section is to characterize the notion of regularity in capacity
by various other properties.
Among those several are functional analytic in nature.
Of special interest is a characterization via relative capacity.
Recall that $\widetilde M$ is the completion of $M$ with respect to
the natural distance and $\partial M = \widetilde M \setminus M$.
The relative capacity is defined on subsets of $\widetilde M$
instead of $M$.
It had been introduced in \cite{AW2}
for an open subset $\Omega$ in $\Ri^d$.
Since it depends on the set $\Omega$ in \cite{AW2} it is called 
relative capacity.
The following definition on manifolds is similar to the Euclidean one.

Let $\mu$ be the trivial extension to $\widetilde M$ of the natural Radon measure
$|\cdot|$ on $M$, that is, for a Borel set $B\subset\widetilde M$ we let $\mu(B)=|B\cap M|$.
For a subset $A\subset\widetilde M$ the {\bf relative capacity} of $A$ (with respect to $M$)
is given by
\[ \rcapp(A) = \{ \|\varphi\|_{H^1(M)}^2 : \varphi\in \widetilde H^1(M) \mbox{ and }
                   \varphi \geq 1\ \mu\mbox{-a.e.\ on a neighbourhood of } A \}
\]
where $\widetilde H^1(M)$ is defined to be the closure of the space
$H^1(M)\cap C_c(\widetilde M)$ in $H^1(M)$.

Note that the relative capacity is the usual capacity as defined in 
\cite{BH} Section~I.8 on the space $\widetilde M$ with respect to the Dirichlet form
$(\psi,\varphi) \mapsto \int_{\widetilde M} \nabla \psi \cdot \nabla \varphi$
and form domain $\widetilde H^1(M)$.
We consider $\widetilde H^1(M)$ instead of $H^1(M)$ in order to
fulfill Condition~(D) in Subsection~I.8.2 of \cite{BH}
and therefore to use the notion of relative quasi-continuity and 
relative quasi-everywhere (r.q.e.).
We do not need that $\widetilde M$ is locally compact, although 
it is a consequence of the embedding theorem of Nash.
In general, however, the completion of a locally compact metric space
is not locally compact.
We are grateful to Robin Nitka for showing us a counter example.

The following characterization of regularity of capacity is our main result in this 
section.
Note that Condition~\ref{tdrum401-4} is formulated completely 
in terms of relative capacity of the boundary $\partial \Omega$.

\begin{thm} \label{tdrum401}
Let $M$ be a connected Riemannian manifold.
Then the following conditions are equivalent.
\begin{tabel}
\item \label{tdrum401-1}
$M$ is regular in capacity.
\item \label{tdrum401-1.5}
The space $C_c^\infty(M)$ is dense in $H^1_0(M) \cap C_0(\widetilde M)$.
\item \label{tdrum401-2}
For every lattice homomorphism $F \colon H^1_0(M) \cap C_0(\widetilde M) \to \Ri$
there exist $c \in \Ri$ and $p \in M$ such that
$F(\varphi) = c \, \varphi(p)$ for all $\varphi \in H^1_0(M) \cap C_0(\widetilde M)$.
\item \label{tdrum401-3}
For every multiplicative functional $\tau$ on the Banach algebra 
$H^1_0(M) \cap C_0(\widetilde M)$ there exists a $p \in M$ such that
$\tau(\varphi) = \varphi(p)$ for all $\varphi \in H^1_0(M) \cap C_0(\widetilde M)$.
\item \label{tdrum401-4}
For every $p\in\partial M$ and $r>0$ one has
$\rcapp(\partial M\cap B_{\widetilde M}(p,r))>0$.
\end{tabel}
\end{thm}

For the proof of Theorem~\ref{tdrum401},  we need  a characterization of
the space $H^1_0(M)$ in terms of the relative capacity.
This result is also of independent interest.

\begin{thm}\label{tdrum520}
Let $M$ be a connected Riemannian manifold.
Then
\begin{equation}
H^1_0(M) = \{ \varphi\in\widetilde H^1(M) : \tilde \varphi=0\mbox{ r.q.e.\ on } \partial M \}
\label{etdrum520;1}
\end{equation}
where $\widetilde \varphi$ denotes the relative quasi-continuous version of $\varphi$.
\end{thm}
\proof\
`$\subset$'.
Since $C^\infty_c(M)\subset H^1(M)\cap C_c(\widetilde M)$
one deduces by closure in $H^1(M)$ that 
$H^1_0(M) \subset \widetilde H^1(M)$.
Let $\varphi\in H^1_0(M)$.
Then it follows from the proof of Proposition~8.2.1 in \cite{BH} that 
there exists a sequence $\varphi_1,\varphi_2,\ldots \in C^\infty_c(M)$ such that
$\lim_{n \to \infty} \varphi_n = \widetilde \varphi$ r.q.e.\ on $\widetilde M$.
So $\widetilde \varphi=0$ r.q.e.\ on $\partial M$.

`$\supset$'.
Let $D^1_0(M)$ denote the right hand side of (\ref{etdrum520;1}).
Let $\varphi\in D^1_0(M)\cap L^\infty(\widetilde M)$.
We may assume that $\varphi \geq 0$.
Then $\varphi \in \widetilde H^1(M)$.
It follows from the definition of $\widetilde H^1(M)$
and the proof of Proposition~8.2.1 in \cite{BH} that there exist
$\varphi_1,\varphi_2,\ldots \in H^1(M) \cap C_c(\widetilde M)$ such that 
$\lim_{n \to \infty} \varphi_n = \tilde \varphi$ in $H^1(M)$ and for all 
$\varepsilon > 0$ there exists an open $U \subset \widetilde M$ such that 
$\rcapp(U) < \varepsilon$ and $\lim_{n \to \infty} \varphi_n = \tilde \varphi$ uniformly 
on $\widetilde M \setminus U$.
We may assume that $0 \leq \varphi_n \leq \|\varphi\|_\infty$
and $\|\varphi_n\|_{H^1(M)} \leq 2 \|\varphi\|_{H^1(M)}$ for all $n \in \Ni$.

Let $\varepsilon \in (0,1]$.
Then there exist $n \in \Ni$ and an open $U \subset \widetilde M$
such that $\| \varphi_n-\varphi\|_{H^1(M)}\leq \varepsilon$, 
$\rcapp(U) < \varepsilon$ and 
$|\varphi_n-\widetilde \varphi|\leq \varepsilon$ uniformly 
on $\widetilde M\setminus U$.
Since $\tilde \varphi = 0$ r.q.e.\ on $\partial M$ there exists an open
$V\subset\widetilde M$ such that 
$ \{ x \in \partial M : \tilde \varphi(x) \neq 0 \} \subset V$ and $\rcapp(V) < \varepsilon$.
Consequently, $\varphi_n\leq \varepsilon$ uniformly on $(\partial M)\setminus W$
where $W=U\cup V$, and $\rcapp(W)\leq \rcapp(U)+\rcapp(V)\leq 2 \varepsilon$.
Let $\chi \in \widetilde H(M)$ be such that $\chi \geq 1$ on $W$ and 
$\|\chi\|_{H^1(M)}^2 < 3 \varepsilon$.
We may assume that
$\chi=1$ pointwise on $W$ and $0\leq\chi\leq 1$ on $M$.
Let $\sigma = (\varphi_n - 2 \varepsilon)^+$ and $\tau = \sigma(\one -\chi)$.
Then $\|\sigma\|_{H^1(M)} \leq 2 \|\varphi\|_{H^1(M)}$ and 
$\|\tau\|_{H^1(M)} \leq 4 \|\varphi\|_{H^1(M)} + 2 \|\varphi\|_\infty$.
Moreover, 
$\|\sigma - \tau\|_2 = \|\sigma \, \chi\|_2 \leq \|\varphi\|_\infty \|\chi\|_2
       \leq 2 \varepsilon^{1/2} \, \|\varphi\|_\infty$.
Then $\supp \tau \subset \supp \sigma \cap W^c$, which is a compact subset of $M$.
So $\tau \in H^1_c(M) \subset H^1_0(M)$.

It follows from the above that 
for all $m \in \Ni$ there exist $\varphi_m,\sigma_m \in H^1(M) \cap L_\infty$
and $\tau_m \in H^1_0(M) \cap L_\infty$ such that 
$\|\varphi - \varphi_m\|_{H^1(M)} \leq \frac{1}{m}$, 
$\|\sigma_m - \tau_m\|_2 \leq \frac{1}{m}$
and $0 \leq \varphi_m - \sigma_m \leq \frac{1}{m}$
for all $m \in \Ni$, and the sequences $\sigma_1,\sigma_2,\ldots$ and 
$\tau_1,\tau_2,\ldots$ are bounded in $H^1(M)$.
We next show that $\tau_1,\tau_2,\ldots$ has a subsequence which converges
to $\varphi$ weakly in $H^1(M)$.

Clearly $\lim \varphi_m = \varphi$ strongly and hence weakly in $H^1(M)$.
The sequences $\sigma_1,\sigma_2,\ldots$ and 
$\tau_1,\tau_2,\ldots$ are bounded in $H^1(M)$.
Hence, by passing to a subsequence if necessary, these sequences are 
weakly convergent in $H^1(M)$.
Since $0 \leq \varphi_m - \sigma_m \leq \frac{1}{m}$
for all $m \in \Ni$ it follows from the 
Lebesgue dominated convergence theorem that 
$\lim \varphi_m - \sigma_m = 0$ in $L_{2,\loc}$
Therefore it follows 
by the uniqueness of the weak limit that 
$\lim \varphi_m - \sigma_m = 0$ weakly in $H^1(M)$.
Because $\lim \sigma_m - \tau_m = 0$ in $L_2$ it follows that 
$\lim \sigma_m - \tau_m = 0$ weakly in $H^1(M)$.
Then $\lim \tau_m = \varphi$ weakly in $H^1(M)$.
So $\varphi \in H^1_0(M)$ and $D^1_0(M) \cap L_\infty(\widetilde M) \subset H^1_0(M)$.

Finally, if $\varphi \in D^1_0(M)$ then 
$(-n) \vee \varphi \wedge n \in D^1_0(M) \cap L_\infty(\widetilde M) \subset H^1_0(M)$
for all $n \in \Ni$ and $\lim (-n) \vee \varphi \wedge n = \varphi$ in $H^1(M)$.
So $\varphi \in H^1_0(M)$.\hfill$\Box$

\ruimte

Finally we prove the characterizations of regular in capacity.

\ruimte

\noindent
{\bf Proof of Theorem~\ref{tdrum401} \hspace{5pt}}\
`\ref{tdrum401-1}$\Rightarrow$\ref{tdrum401-1.5}'.
Let $\varphi \in H^1_0(M) \cap C_0(\widetilde M)$ and $\varepsilon>0$.
We may assume that $\varphi \geq 0$.
Since $\varphi \in C_0(\widetilde M)$ there exits a compact 
$K\subset\widetilde M$ such that $\varphi(q)<\varepsilon$ for all
$q \in \widetilde M \setminus K$.
Moreover, $\varphi(q) = 0$ for all $q\in \partial M$ since 
$M$ is regular in capacity.
Let $U = \{ q \in \widetilde M : \varphi(q) < \varepsilon \} $.
Then $U$ is open and $\widetilde M \setminus M \subset U$.
Moreover, $\supp (\varphi - \varepsilon)^+ \subset (\widetilde M \setminus U) \cap K$
and hence compact.
But $(\widetilde M \setminus U) \cap K \subset M$.
Therefore $(\varphi - \varepsilon)^+ \in H^1_0(M) \cap C_c(M)$.
Using a partition of the unity one deduces that 
$C_c^\infty(M)$ is dense in $H^1_0(M) \cap C_c(M)$.
Finally, $\lim_{\varepsilon \downarrow 0} (\varphi - \varepsilon)^+ = \varphi$
in $H^1_0(M) \cap C_0(\widetilde M)$.
So $C_c^\infty(M)$ is dense in $H^1_0(M) \cap C_0(\widetilde M)$.

`\ref{tdrum401-1.5}$\Rightarrow$\ref{tdrum401-1}'.
Suppose $M$ is not regular in capacity.
Then there are $\varphi \in H^1_0(M) \cap C_0(\widetilde M)$ and $p \in \partial M$
such that $\varphi(p) \neq 0$.
Then $\|\varphi - \psi\|_{C_0(\widetilde M)} \geq |\varphi(p)|$ for all 
$\psi \in C_c^\infty(M)$, so $C_c^\infty(M)$ is not dense in $H^1_0(M) \cap C_0(\widetilde M)$.

`\ref{tdrum401-1.5}$\Rightarrow$\ref{tdrum401-2}'.
Let $F \colon H^1_0(M) \cap C_0(\widetilde M) \to \Ri$ be a lattice 
homomorphism.
Then $F$ is continuous by \cite{Schae} Theorem~V.5.5(ii).
Arguing as at the end of the proof of Lemma~\ref{ldrum201} 
it follows from Lemma~\ref{ldrum201.3} that 
there are $c \in \Ri$ and $p \in M$ such that 
$F(\varphi) = c \, \varphi(p)$ for all $\varphi \in C_c^\infty(M)$.
Since $F$ is continuous and 
$C_c^\infty(M)$ is dense in $H^1_0(M) \cap C_0(\widetilde M)$ it 
follows that $F(\varphi) = c \, \varphi(p)$ for all 
$\varphi \in H^1_0(M) \cap C_0(\widetilde M)$.

`\ref{tdrum401-2}$\Rightarrow$\ref{tdrum401-1}'.
Suppose $M$ is not regular in capacity.
Then there are $\psi \in H^1_0(M) \cap C_0(\widetilde M)$ and 
$p \in \partial M$ such that $\psi(p) \neq 0$.
Define $F \colon H^1_0(M) \cap C_0(\widetilde M) \to \Ri$ by
$F(\varphi) = \varphi(p)$.
Then $F$ is a continuous lattice homomorphism.
So by assumption there are $q \in M$ and $c \in \Ri$ such that $F(\varphi) = c \, \varphi(q)$
for all $\varphi \in H^1_0(M) \cap C_0(\widetilde M)$.
Let $\chi \in C_c^\infty(M)$ such that $\chi(q) = 1$.
Then $\psi(\one - \chi) = \psi-\psi\chi \in H^1_0(M) \cap C_0(\widetilde M)$.
Therefore 
\[
0 \neq \psi(p) 
= (\psi(\one - \chi))(p) 
= F(\psi(\one - \chi))
= c \, (\psi(\one - \chi))(q)
= 0
\;\;\; .\]
This is a contradiction.

`\ref{tdrum401-1.5}$\Rightarrow$\ref{tdrum401-3}'.
Let $\tau \colon H^1_0(M) \cap C_0(\widetilde M) \to \Ci$ be a (non-zero) multiplicative functional.
Then $\tau$ is continuous by \cite{HR1}, Theorem C.21.
Therefore it follows by Condition~\ref{tdrum401-1.5} that 
$\tau|_{C_c^\infty(M)} \colon C_c^\infty(M) \to \Ci$ is a 
(non-zero) multiplicative functional.
Let $p,q \in \supp \tau|_{C_c^\infty(M)}$ with $p \neq q$ and let
$U$ and $V$ be two disjoint open neighbourhoods of $p$ and $q$ respectively.
Then there exist $\varphi \in C_c^\infty(U)$ and $\psi \in C_c^\infty(V)$ such that 
$\tau(\varphi) \neq 0$ and $\tau(\psi) \neq 0$.
But then $\varphi \, \psi = 0$ and 
\[
0 = \tau(\varphi \, \psi) 
= \tau(\varphi) \, \tau(\psi)
\neq 0
\;\;\; .\]
This is a contradiction.
So there exists a $p \in M$ such that $\supp \tau|_{C_c^\infty(M)} = \{ p \} $.

Next we show that $\tau$ is positive.
Let $\varphi \in C_c^\infty(M)$ and suppose that $\varphi \geq 0$.
If $\varphi(p) > 0$ then there exist $\psi \in C_c^\infty(M)$ and a
neighbourhood $V$ of $p$ such that $\varphi|_V = \psi^2|_V$.
Then $\tau(\varphi) = \tau(\psi^2) = \tau(\psi)^2 \geq 0$.
Alternatively, suppose that $\varphi(p) = 0$.
Let $V$ be a relative compact neighbourhood of $p$.
Then by continuity there exists a $c > 0$ such that 
$|\tau(\psi)| \leq c \, \|\psi\|_{W^{1,\infty}(V)}$
for all $\psi \in C_c^\infty(V)$.
We may assume that $\supp \varphi \subset V$.
Then $\lim_{\varepsilon \downarrow 0} (\varphi - \varepsilon)^+ = \varphi$ in 
$W^{1,\infty}(V)$, so by regularizing it follows that there are 
$\varphi_1,\varphi_2,\ldots \in C_c^\infty(V)$ such that $\lim \varphi_k = \varphi$
in $W^{1,\infty}(V)$ and $\varphi_k$ vanishes in a neighbourhood of $p$ for all
$k \in \Ni$.
Then $\tau(\varphi) = \lim \tau(\varphi_k) = 0$.

Now it follows from Lemma~\ref{ldrum201.3} that there are $c \in [0,\infty)$ and $p \in M$
such that $\tau(\varphi) = c \, \tau(\varphi)$ for all $\varphi \in C_c^\infty(M)$.
Then $c^2 = 1$ and since $\tau \neq 0$ it follows that $c = 1$.
Since $C_c^\infty(M)$ is dense in $H^1_0(M) \cap C_0(\widetilde M)$ one establishes
that $\tau(\varphi) = \varphi(p)$ for all $\varphi \in H^1_0(M) \cap C_0(\widetilde M)$.

`\ref{tdrum401-3}$\Rightarrow$\ref{tdrum401-1}'.
This proof is similar to the proof \ref{tdrum401-2}$\Rightarrow$\ref{tdrum401-1}.

`\ref{tdrum401-1}$\Rightarrow$\ref{tdrum401-4}'.
Assume that there exist $p\in \partial M$ and $r>0$ such that
$\rcapp(B_{\widetilde M}(p,r)\cap \partial M) = 0$.
Then there exist a $\widetilde M$-open neighbourhood
$U$ of $B_{\widetilde M}(p,r)\cap \partial M$ and a function $\chi\in\widetilde H^1(M)$ such that
$\chi\geq 1$ a.e.\ on $U\cap M$.
Let $\rho\in (0,r)$ be such that $M\setminus B_{\widetilde M}(p,\rho)\ne\emptyset$.
Define  $\psi \colon \widetilde M\to\Ri$  by 
$\psi(q) = d_{\widetilde M}(q,\widetilde M\setminus B(p\,;\rho))$.
Then $\psi \in C(\widetilde M)$ and $\psi|_M \in W^{1,\infty}(M)$.
Set $\varphi=\chi \, \psi$.
Then $\varphi\in C(\widetilde M)$ and by an elementary argument
one deduces that $\varphi|_M\in\widetilde H^1(M)$.
Moreover, $\varphi = 0$ r.q.e.\ on $\partial M$.
By Theorem \ref{tdrum520} it follows that $\varphi|_M\in H^1_0(M)$.
So $\varphi\in H^1_0(M) \cap C_0(\widetilde M)$.
Since $\varphi(p)\ne 0$ it follows from the definition that the
manifold $M$ is not regular in capacity.

`\ref{tdrum401-4}$\Rightarrow$\ref{tdrum401-1}'.
If $M$ is not regular in capacity then there exist $\varphi\in H^1_0(M) \cap C_0(\widetilde M)$ and
$p\in\partial M$ such that $\varphi(p) \neq 0$.
Without loss of
generality we may assume that $\varphi(p)=2$.
Let $r\in (0,1)$ be such that $\varphi\geq 1$ on $B_{\widetilde M}(p,r)$.
Since $\varphi \in H^1_0(M)$ one deduces from 
Theorem~\ref{tdrum520}  that $\varphi=0$ r.q.e.\ on $\partial M$.
Then 
$\rcapp( B_{\widetilde M}(p,r)\cap \partial M )
\leq \rcapp( \{ q \in \partial M : \varphi(q) \neq 0 \} ) = 0$.\hfill$\Box$


\bibliography{refbib}

\end{document}